\documentclass[]{elsarticle}
\usepackage{caption}
\usepackage{subcaption}
\usepackage{hyperref}
\usepackage{amsmath}
\usepackage{amssymb}
\usepackage{mathtools}
\usepackage{cleveref}
\usepackage{bm}
\usepackage{amsbsy}
\usepackage{multirow}

\usepackage{a4wide}
\usepackage{color}

\newcommand{\xvec}[1]{\pmb{#1}}

\biboptions{sort&compress}

\begin{document}
\begin{frontmatter}
\title{A data-based reduced-order model for dynamic simulation and control of district-heating networks}

\author{Mengting Jiang\corref{correspondingauthor}}\cortext[correspondingauthor]{Corresponding author}\ead{m.jiang1@tue.nl}
\author{Michel Speetjens}
\author{Camilo Rindt}
\author{David Smeulders}
\address{Eindhoven University of Technology, Mechanical Engineering Department}

\begin{abstract}

This study concerns the development of a data-based compact model for the prediction of the fluid temperature evolution in district heating (DH) pipeline networks. This so-called ``reduced-order model'' (ROM) is obtained from reduction of the conservation law for energy for each pipe segment to a semi-analytical input-output relation between the pipe outlet temperature and the pipe inlet and ground temperatures that can be identified from training data. The ROM basically is valid for generic pipe configurations involving 3D unsteady heat transfer and 3D steady flow as long as heat-transfer mechanisms are linearly dependent on the temperature field. Moreover, the training data can be generated by physics-based computational ``full-order'' models (FOMs) yet also by (calibration) experiments or field measurements.
Performance tests using computational training data for single 1D pipe configuration demonstrate that the ROM (i) can be successfully identified and (ii) can accurately describe the response of the
outlet temperature to arbitrary input profiles for inlet and ground temperatures. Application of the ROM to two case studies, i.e. fast simulation of a small DH network and design of a
controller for user-defined temperature regulation of a DH system, demonstrate its predictive ability and efficiency also for realistic systems. Dedicated cost analyses further reveal that the
ROM may significantly reduce the computational costs compared to FOMs by (up to) orders of magnitude for higher-dimensional pipe configurations. These findings advance the proposed ROM as a robust and efficient simulation tool for practical DH systems with a far greater predictive ability than existing compact models.

\end{abstract}

\end{frontmatter}




\section{Introduction}


Studies from the International Energy Agency (IEA) show that domestic energy demand mainly concerns heating of buildings and, since
this is primarily satisfied by fossil fuels, contributes by more than 40\% to the global energy-related carbon dioxide $(CO_2)$ emissions \citep{Ferroukhi2020}. This energy demand is predicted to only increase in the future \citep{Capuano2018}. Thus developing more sustainable
ways for heating of the built environment is critical to mitigate environmental issues such as e.g. global warming and pollution caused by using fossil fuels.

District-heating (DH) systems distributing thermal energy extracted from geothermal reservoirs via water flows in
pipeline networks is a promising solution for the sustainable heating and cooling of buildings. The study by \citep{Connolly2014} in fact shows
that DH systems can achieve the same reduction in primary energy demand and associated $CO_2$ emissions at a lower (economical) cost
compared with other alternatives such as electrification. However, the market share of DH systems has in most countries nonetheless remained low mainly due to its inferior level of thermal comfort compared to other technologies \citep{Talebi2016}.
This is largely caused by (i) the centralized and inflexible nature of conventional DH systems due to their reliance on a single energy source with fixed supply temperatures and (ii) the only limited control users
have on the heat supply. Future DH systems should offer much greater thermal comfort by facilitating dynamic energy demands and supply temperatures. Future DH systems must therefore enable active interaction between users and the network by making them, besides energy {\it consumers}, also energy {\it producers}, i.e. so-called ``prosumers'' \citep{Lund2014}. Future DH systems must furthermore allow
for decentralized energy supply from multiple and intermittent energy sources (including energy storage) as well as integration in combined heat and power networks \citep{Lund2014,Li2015a,Li2017}. Crucial to these ends is flexible and pro-active operation and control of DH systems based on feedback and input from consumers \citep{Vandermeulen2018}.

A first step towards pro-active operation and control of DH systems exists in so-called ``optimal scheduling'', that is, determining an optimal interaction scheme between the
various components in the (integrated) energy network on the basis of anticipated energy supplies and demands over a given period of time \cite{Gu2017,Deng2017,Wang2019,Merkert2020a}.
This is in essence {\it open-loop control} by running the DH network according to an {\it a priori} determined procedure without adjustment or intervention during its
actual operation \cite{ControlTheoryEng}. Thus optimal scheduling is incapable of responding to unforeseen behaviour and disturbances.
%
%
True pro-active operation of the DH network, accounting for both expected and unforeseen events, requires {\it closed-loop control} by regulating the system components and their interactions on the basis
of feedback on the actual system behaviour during operation. Suitable control strategies for this purpose -- finding widespread and successful application in the intimately related field of process
control in industry -- include the conventional PID controllers and advanced methods as Model Predictive Control (MPC) \citep{ControlTheoryEng,Camacho2013,Corriou2018}.

Essential for both open-loop and closed-loop control of DH networks is a computational model that is compact, efficient and fast and at the same time sufficiently accurate to capture the
relevant dynamics. (A key factor in the latter is thermal inertia of the DH network stemming (mainly) from the pipelines, buildings and working fluid \cite{Gu2017,Wang2019}.)
Various optimization algorithms for optimal scheduling are e.g. readily available in the literature \citep{Talebi2016,Gu2017,Deng2017,Wang2019,Merkert2020a}; efficient and adequate models are typically the main challenge \citep{Li2017}. Such models are also critical for proper design and fine-tuning of e.g. PID controllers (``controller synthesis'') \citep{ControlTheoryEng}.
Real-time control by MPC, arguably, is most critically dependent on an efficient model, since this control strategy relies on repeated predictions of future behaviour from feedback on the intermediate state of the system and subsequent determination of a control action at much faster times scales than the system dynamics \citep{Camacho2013,Corriou2018}.

The most detailed and accurate computational models for DH networks result from spatial discretization of the conservation laws for flow and heat transfer in the system components
by e.g. finite-element or finite-volume methods \cite{Hirsch2006}. However, such models are computationally expensive and thus typically ill-suited for the above control purposes
and basically disqualify for advanced closed-loop control (by e.g. MPC) \cite{Palsson2000, Gabrielaitiene2008, Sartor2015, Stevanovic2009}. This motivated the development of a range
of compact models for DH systems and then mostly pipe models for the integration in network models.

The ``node method'' originally by \cite{Benonysson1991} is the most common approach and estimates the pipe outlet temperature from the inlet temperature by accounting
for the travel time of fluid parcels within the pipe as well as thermal inertia/loss of/via the pipe wall \cite{Palsson2000,Gabrielaitiene2008,zheng2017}.
The ``characteristic method'' by \cite{Stevanovic2009} considers the actual unsteady 1D energy balance relative to the flow and is more accurate than the node method. However, the necessary discretization of the flow region renders this method computationally more expensive.
The ``function method'' by \cite{zheng2017} seeks to unite the advantages of node and characteristic methods by relating inlet and outlet temperatures through an analytical function deriving from coupled unsteady 1D energy balances for fluid and pipe wall. The aggregated model by \citep{larsen2002aggregated} is a compact network model that reduces the DH system to an equivalent simpler network of pipe segments each described by an integral quasi-steady energy balance including finite travel time and thermal loss (yet without thermal inertia).

The principal shortcoming of existing compact models as those reviewed above is their reliance on (oversimplified) energy balances and the resulting limited predictive ability. The main objective of the present study is development of a compact pipe model that offers far greater physical validity at comparable computational cost. This to a certain extent expands on the function method by \cite{zheng2017} by generalizing the analytical relation between inlet and outlet temperatures to more complex pipe configurations involving 3D unsteady heat transfer and 3D steady flow. The proposed model is to be founded on input-output relations for generic linear dynamical systems and thus basically the only condition for its validity is heat-transfer mechanisms being linearly dependent on the temperature field. Moreover, the proposed model can be identified entirely from data generated by physics-based computational models yet also by (calibration) experiments or field measurements.

The paper is organized as follows. \Cref{numericalmethods} defines the pipe configuration and the corresponding computational model, denoted ``full-order model'' (FOM) hereafter, to simulate its thermal behaviour. The compact model, denoted ``reduced-order model'' (ROM) hereafter, and its identification from so-called ``training data'' (generated by the FOM) is discussed in \Cref{ROM}. Performance and validation of the ROM are investigated in \Cref{discussion}. Application of the ROM for the simulation and control of DH systems is demonstrated in \Cref{application} including an analysis of computational costs. Conclusions and future work are discussed in \Cref{conclusion}.

\section{Problem definition and full-order model (FOM)}\label{numericalmethods}

\subsection{System configuration and governing equations}

\begin{figure}[htbp!]
\begin{centering}
	\begin{subfigure}[b]{0.38\textwidth}
\centering
   \includegraphics[width=\textwidth]{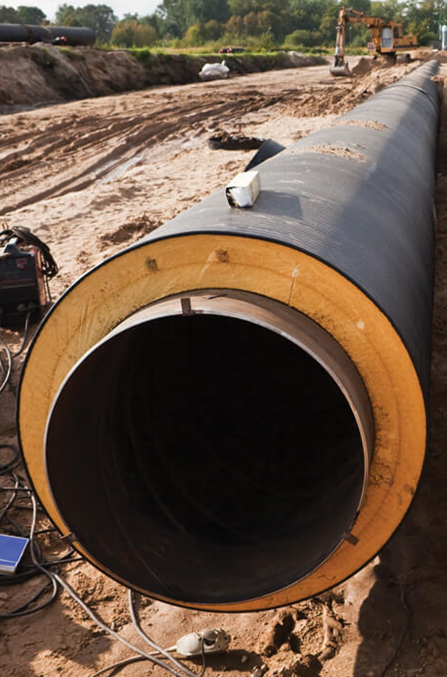}
	\caption{}
	\label{picture-pipe}
	\end{subfigure}
	\hfill
	\begin{subfigure}[b]{0.6\textwidth}
\centering
   \includegraphics[width=0.8\textwidth]{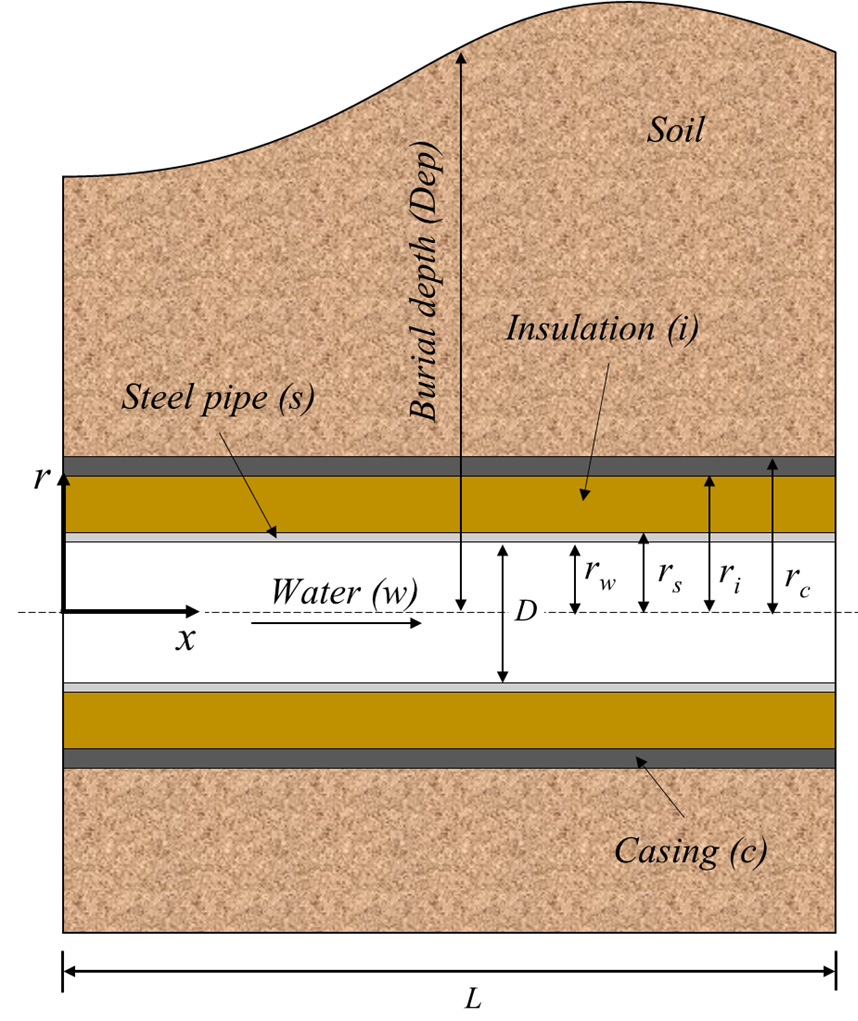}
	\caption{}
	\label{sketch-pipe}
	\end{subfigure}
\caption{Typical underground insulated DH pipe: (a) pipe before installation; (b) schematic of installed configuration. Panel (a): \url{https://www.bcc-extrusion.com/applications/pipe/po-pipes/district-heating-pipes.html}.}
\end{centering}

\end{figure}

The problem of interest is the transfer and distribution of heat via the flow of water through DH networks of underground insulated pipes as shown in \Cref{picture-pipe}. Each pipe segment of length $L$ has a configuration following the vertical cross-section in \Cref{sketch-pipe} and consists of a steel pipe with inner and outer radius $r_w$ and $r_s$, respectively, covered by a cylindrical insulation of outer radius $r_i$ and a cylindrical casing of outer radius $r_c$ buried in the soil at depth $Dep$.
Modeling of the thermo-hydraulic behavior in such a pipe segment is done by the model proposed by Ref.~\cite{Palsson2000} and is based on the following assumptions:
\begin{itemize}
\item the system is in a quasi-dynamic state, i.e. the flow and heat transfer are steady and unsteady, respectively;
\item the fluid (i.e. water) is incompressible and has constant thermo-physical properties;
\item the flow is uniform in both radial ($r$) and axial ($x$) directions;
\item the heat transfer in flow and solid regions is quasi-2D, i.e. the temperature field exhibits full spatial variation in axial ($x$) direction and is component-wise uniform in radial ($r$) direction.
\end{itemize}
%
%
These conditions reduce the fluid motion to a uniform axial flow of magnitude $u$. The thermal behavior in the pipe configuration in \Cref{sketch-pipe} is described in terms of the temperature as a function of axial coordinate $x$ and time $t$, i.e. $T=T(x,t)$, which is governed by the following energy balances subject to the above assumptions:
\begin{equation}\label{eqwater}
\frac{\partial T^w}{\partial t} + u \frac{\partial T^w}{\partial x} = \alpha^w \frac{\partial^2 T^w}{\partial x^2} + \frac{1}{R_{ws} V_w \rho_w C_{p_w}}(T^s - T^w),
\end{equation}
\begin{equation}\label{eqsteel}
\frac{\partial T^s}{\partial t} = \alpha^s \frac{\partial^2 T^s}{\partial x^2} + \frac{1}{R_{ws} V_s \rho_s C_{p_s}}(T^w - T^s) - \frac{1}{R_{si} V_s \rho_s C_{p_s}}(T^s - T^i),
\end{equation}
\begin{equation}\label{eqinsu}
\frac{\partial T^i}{\partial t} = \alpha^i \frac{\partial^2 T^i}{\partial x^2} + \frac{1}{R_{si} V_i \rho_i C_{p_i}}(T^s - T^i) - \frac{1}{R_{ic} V_i \rho_i C_{p_i}}(T^i - T^c),
\end{equation}
\begin{equation}\label{eqcasing}
\frac{\partial T^c}{\partial t} = \alpha^c \frac{\partial^2 T^c}{\partial x^2} + \frac{1}{R_{ic} V_c \rho_c C_{p_c}}(T^i - T^c) - \frac{1}{R_{cg} V_c \rho_c C_{p_c}}(T^c - T_g),
\end{equation}
where superscripts for temperature $T$ refer to the different subregions of the system, i.e. water (``w''), steel pipe (``s''), insulation (``i'') and casing (``c''), and $T_g$ corresponds with the uniform ground temperature. The set of PDEs \eqref{eqwater}-\eqref{eqcasing} is subject to adiabatic boundary conditions for $T^{w,s,i,c}$, a Dirichlet boundary condition for $T_w$ and uniform initial conditions following
\begin{equation}
\frac{\partial T^{s,i,c}(x,t)}{\partial x}|_{x=0},\quad\quad T^w(0,t) = T_{in}(t),\quad\quad T^{w,s,i,c}(x,0) = T_0,
\label{BCs}
\end{equation}
with $T_{in}$ the (variable) inlet temperature of the fluid.

The trailing term on the LHS of \eqref{eqwater} describes axial heat transfer by convection due to the uniform flow $u$ and the leading terms on the RHS of \eqref{eqwater}-\eqref{eqcasing} describes axial heat transfer by conduction parameterized by thermal diffusivity $\alpha$ of the respective subregions. The remaining terms on the RHS of \eqref{eqwater}-\eqref{eqcasing} describe the radial heat transfer between the subregions, where $(V,\rho,C_p)$ are the volume, density and heat capacity of each element, respectively, and $R$ define thermal resistances between adjacent regions according to the subscript (e.g. $R_{ws}$ for the heat exchange between water and steel pipe). The thermal resistances are given by:
\begin{eqnarray}
R_{ws}=\frac{1}{2\pi r_wL h_w}+\frac{ln\left((r_w+r_{ws})/r_w\right)}{2\pi L k_s},\quad\quad
R_{si}=\frac{ln\left(r_s/(r_s-r_{ws})\right)}{2\pi L k_s}+\frac{ln\left((r_s+r_{si})/r_s\right)}{2\pi L k_i},\\\nonumber\\
R_{ic}=\frac{ln\left(r_i/(r_i-r_{si})\right)}{2\pi L k_i}+\frac{ln\left((r_i+r_{ic})/r_i\right)}{2\pi L k_c},\quad\quad
R_{cg}=\frac{ln\left(r_c/(r_c-r_{ic})\right)}{2\pi L k_c}+R_{soil},
\end{eqnarray}
with $r_{ws}=(r_s-r_w)/2$, $r_{si}=(r_i-r_s)/2$ and $r_{ic}=(r_c-r_i)/2$, where $h_w$ is the convective heat-transfer coefficient for the flow, $k$ is the thermal conductivity of each element and
$R_{soil}$ is the thermal resistance of the soil. Coefficient $h_w$ is computed from the Nusselt relation \citep{Bergman2011,Clamond2009}
\begin{eqnarray}
Nu=\frac{h_wD}{k_w}= \frac{\left(f/8\right) \left(Re-1000\right)P_r}{1+12.7\sqrt{(f/8)\left(Pr^\frac{2}{3}-1\right)}},\quad\quad \frac{1}{\sqrt f}=-2 log_{10}\left(\frac{\left(\varepsilon/D\right)}{3.7}+\frac{2.51}{Re}\frac{1}{\sqrt f}\right),
\label{Nusselt}
\end{eqnarray}
describing the Nusselt number $Nu$ as a function of the Darcy-Weisbach friction factor $f$, defined implicitly by the right relation in \eqref{Nusselt}, the Reynolds
number $Re=\rho_w u D/\mu_w$ and the Prandtl number $Pr=C_{p_w}\mu_w/k_w$, with $\varepsilon$ the pipe roughness, $D$ the pipe diameter and $\mu_w$ the dynamic viscosity of the fluid.
The thermal resistance of the soil is evaluated via relation \citep{Jianguang2018}
\begin{equation}
R_{soil}=\frac{ln\left(Dep/r_c+\sqrt{\left(Dep/r_c\right)^2-1}\right)}{2\pi L k_s},
\end{equation}
with $Dep$ the burial depth of the pipe (relative to the centerline).

\subsection{Full-order model (FOM)}\label{FOM}

\begin{figure}[htbp!]
\centering
\includegraphics[width=12cm, height=5.5cm]{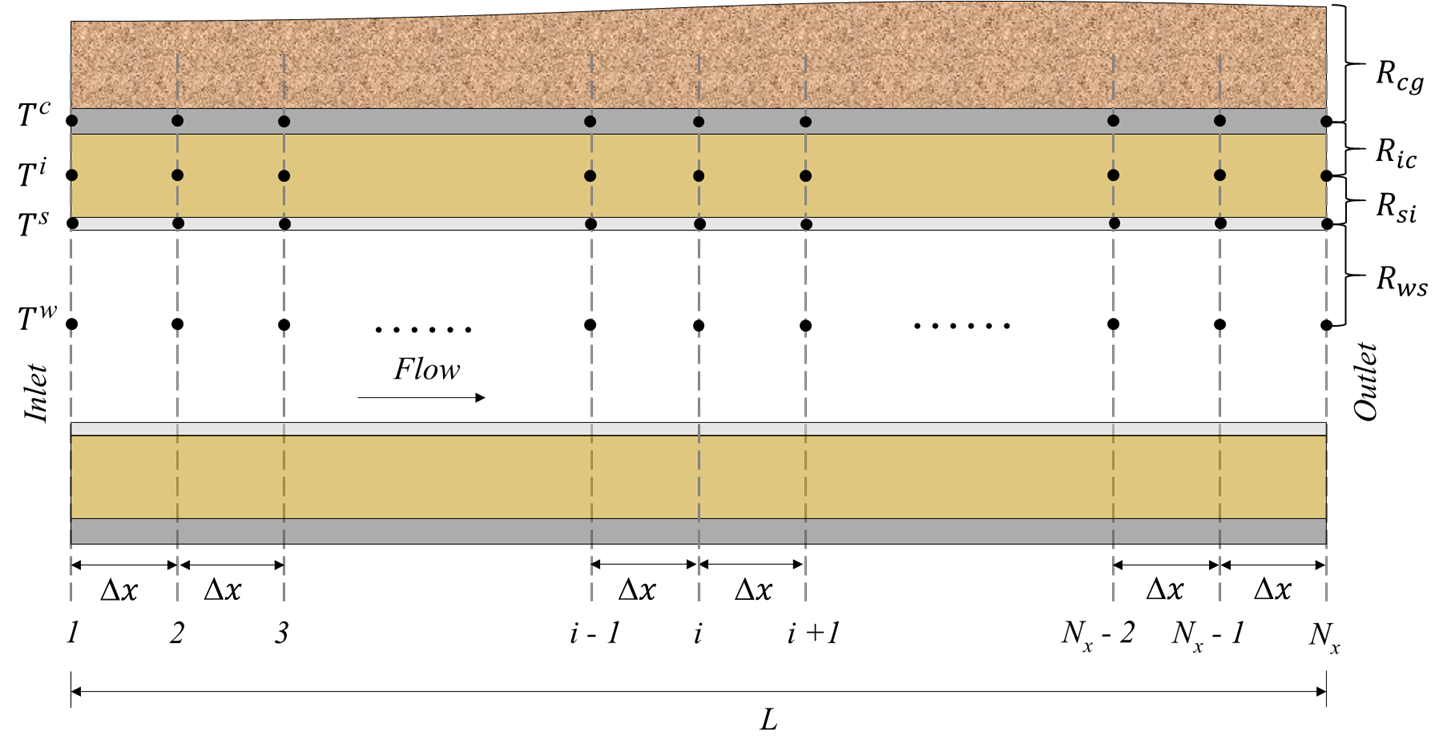}
\caption{Spatial discretization of the pipe segment shown in \Cref{sketch-pipe} by partitioning into $N_x$ cells of length $\Delta x$ in the stream-wise ($x$) direction.}
\label{pipe-discretization}
\end{figure}

The FOM is built from spatio-temporal discretization of PDEs (\ref{eqwater})-(\ref{eqcasing}) using conventional finite-difference methods for space and time \citep{Jayanti2018,Thomas2013}. To this end the pipe is partitioned into cells of length $\Delta x = L/(N_x-1)$ and separated by nodes $x_i=(i-1)\Delta x$, with $1\leq i \leq N_x$ and $N_x\geq 2$ the number of nodes, following \Cref{pipe-discretization} for discretization in space with a first-order upwind scheme for the convection term and a central-difference scheme for the conduction terms. The corresponding time evolution is discretized with an explicit first-order Euler scheme using a finite time step $\Delta t$ and discrete time levels $t_j = j\Delta t$, with $j\geq 0$. This results in a FOM that propagates the temperature of the subregions at position $x_i$ from
time level $t_j$ to time level $t_{j+1}$ via
\begin{eqnarray}
T_{i,j+1}^w &=& T_{i,j}^w + \Delta t\left[\alpha^w \frac{T_{i+1,j}^w - 2T_{i,j}^w + T_{i-1,j}^w}{\Delta x^2} - u_j\frac{T_{i,j}^w - T_{i-1,j}^w}{\Delta x}- \frac{\Delta x(T_{i,j}^w - T_{i,j}^s)}{R_{ws} V_w \rho_w C_{p_w}L}\right],
\label{discretizewater}\\\nonumber\\
T_{i,j+1}^s &=& T_{i,j}^s + \Delta t\left[\alpha^s \frac{T_{i+1,j}^s - 2T_{i,j}^s + T_{i-1,j}^s}{\Delta x^2} + \frac{\Delta x(T_{i,j}^w - T_{i,j}^s)}{R_{ws} V_s \rho_s C_{p_s}L}
- \frac{\Delta x(T_{i,j}^s - T_{i,j}^i)}{R_{si} V_s \rho_s C_{p_s}L}\right],\label{discretizepipe}\\\nonumber\\
T_{i,j+1}^i &=& T_{i,j}^i + \Delta t\left[\alpha^i \frac{T_{i+1,j}^i - 2T_{i,j}^i + T_{i-1,j}^i}{\Delta x^2}
+ \frac{\Delta x(T_{i,j}^s - T_{i,j}^i)}{R_{si} V_i \rho_i C_{p_i}L} - \frac{\Delta x(T_{i,j}^i - T_{i,j}^c)}{R_{ic} V_i \rho_i C_{p_i}L}\right],\label{discretizeinsulation}\\\nonumber\\
T_{i,j+1}^c &=& T_{i,j}^c + \Delta t\left[\alpha^c \frac{T_{i+1,j}^c - 2T_{i,j}^c + T_{i-1,j}^c}{\Delta x^2}
+ \frac{\Delta x(T_{i,j}^i - T_{i,j}^c)}{R_{ic} V_c \rho_c C_{p_c}L} - \frac{\Delta x(T_{i,j}^c - T_g)}{R_{cg} V_c \rho_c C_{p_c}L}\right],\label{discretizecasing}
\end{eqnarray}
subject to the discrete counterparts of conditions \eqref{BCs}, i.e.
\begin{equation}
T^{s,i,c}_{2,j}-T^{s,i,c}_{1,j} = 0,\quad\quad
T^w_{1,j} = T_in(t_j), \quad\quad
T_{i,j}^{w,s,i,c} = T_0,
\label{BCs2}
\end{equation}
where $T_{i,j}^{w,s,i,c} = T^{w,s,i,c}(x_i,t_j)$.

An important computational advantage of the explicit time discretisation is that this decouples the temperature evolutions by enabling
evaluation of the unknown temperature $T^{w,s,i,c}_{i,j+1}$ in each position $x_i$ and subregion from the individual relations
\eqref{discretizewater}-\eqref{discretizecasing} using the fully known temperature distribution $T^{w,s,i,c}$ at time level $t_{j}$. This
admits a very efficient propagation in time even for variable flow and/or thermo-physical properties. Disadvantage of this approach is a
conditional numerical stability, resulting in the following stability criterion for the time step \citep{Thomas2013}:
\begin{equation}\label{stability}
\Delta t \leq \Delta t_{max},\quad\quad\Delta t_{max} = \frac{1}{u/\Delta x + 2\alpha_{min}/\Delta x ^2}.
\end{equation}
with $\alpha_{min} = \min \alpha^{w,s,i,c}$.

The FOM \eqref{discretizewater}-\eqref{discretizecasing} can be collected into a matrix-vector relation for the entire system, i.e.
\begin{eqnarray}
\xvec{T}_{j+1} = \xvec{T}_{j} + \Delta t\left( \xvec{A}\xvec{T}_{j} + T_{in}\xvec{b}_1 + T_{g}\xvec{b}_2\right),
\label{GlobalFOM}
\end{eqnarray}
with column vector
\begin{eqnarray}
\xvec{T}_j = \left[T^w_{1,j},\dots, T^w_{N_x,j}\quad\cdots\quad T^i(x_{N_x},t_j),\dots,T^c(x_{N_x},t_j)\right]^\dagger,
\label{StateVector}
\end{eqnarray}
representing the global state (i.e. the temperatures of all subregions in all nodes) at time level $t_j$ (symbol $^\dagger$ indicates
transpose). Here system matrix $\xvec{A}$ captures the spatially-dependent terms in \eqref{discretizewater}-\eqref{discretizecasing} and system vectors $\xvec{b}_{1,2}$
capture the impact of the inlet condition $T_{in}$ and the soil temperature $T_g$ and the global temperature evolution.
The global FOM \eqref{GlobalFOM} paves the way to the ROM to be developed hereafter due to the fact that it results from temporal discretization of
\begin{equation}\label{semi-discrete}
\frac{d\xvec{T}}{dt} = \xvec{A}\xvec{T} + \xvec{B}\xvec{g},\quad\quad
\xvec{B} = [\xvec{b}_1\;\xvec{b}_2],\quad\quad
\xvec{g} = [T_{in}\;T_g]^\dagger,
\end{equation}
with the beforementioned explicit Euler scheme.
Semi-discrete model \eqref{semi-discrete} namely describes the continuous evolution of state vector \eqref{StateVector} in time $t$, i.e. $\xvec{T}=\xvec{T}(t)$, and constitutes a
generic ``linear time-invariant (LTI) system'' from control theory, with $\xvec{A}$ and $\xvec{B}$ the system and input matrices, respectively, and $\xvec{g}$ the input \cite{ControlTheoryEng}.
Thus \eqref{semi-discrete} bridges the gap to methods and concepts from control theory for the development of the ROM. This is elaborated in \Cref{ROM}.

Spatial discretisation of the (set of) energy balance(s) of more complex pipe configurations involving 3D unsteady heat transfer and 3D steady flow (upon representation of the inlet conditions by a mean $T_{in}$) yields a semi-discrete model with an identical structure as \eqref{semi-discrete} as long as heat-transfer mechanisms are linearly dependent on the temperature field \cite{Hirsch2006}. This has the important implication that the ROM to be developed in \Cref{ROM}, though here demonstrated for the 1D pipe configuration according to \Cref{numericalmethods}, is in fact valid for a much wider range of systems.

\section{Reduced-order model (ROM)}\label{ROM}

The basis for the ROM is an input-output relation associated with the FOM that directly expresses the fluid temperature at the pipe outlet $x=L$,
i.e. $T_{out}(t) = T^w(L,t)$, as a function of the input $\xvec{g}$ of the LTI system \eqref{semi-discrete}. Such input-output relations
are a common way to model practical systems and processes, because, first, they admit efficient description of the behaviour particularly for linear systems and, second, they often admit identification from data \cite{ControlTheoryEng}. Both features will be exploited here as well. Reduction of the FOM to an input-output relation is explained in \Cref{Reduction}. Construction of a ROM for generic input $\xvec{g}$ and its identification from data is elaborated in \Cref{UnitResponse} and \Cref{FunctionIdentify}, respectively.

\subsection{Reduction of the FOM to an input-output relation}\label{Reduction}

\begin{figure}[h!]
\centering
\includegraphics[width=8cm, height=3.5cm]{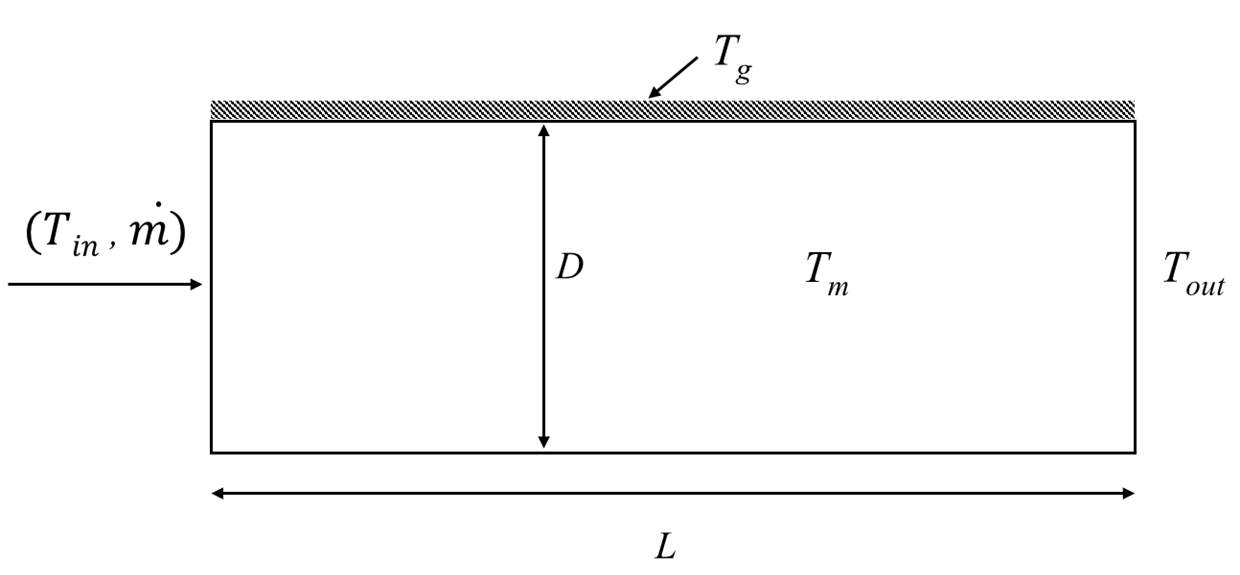}
\caption{Simplified pipe segment with length $L$ and diameter $D$ for illustration of the concept of an input-output relation between inlet ($T_{in}$) and outlet ($T_{out}$) temperatures.}
\label{pipe-simplified}
\end{figure}

Consider for illustration of the concept of an input-output relation the simplified system depicted in \Cref{pipe-simplified}. The integral energy balance of this simplified problem is:
\begin{equation}\label{eqsimplified}
m C_{p_w}\frac{dT_m}{dt}=\dot{m}C_{p_w}(T_{in}-T_{out}) + h_w\mathcal{A}(T_g-T_m),
\end{equation}
with $T_m$ the mean fluid temperature, $m$ the mass content of fluid inside the pipe, $\dot{m}$ the mass flux, $T_{in}$ the temperature at the pipe inlet, $T_{out}$ the temperature at the
pipe outlet, $T_g$ the ground (and wall) temperature, $h_w$ the heat-transfer coefficient from pipe wall (area $\mathcal{A} = \pi D L$) to bulk flow.

Assuming $T_{out} = T_m$ as relation between outlet and mean temperature -- which, to good approximation, holds for sufficiently high thermal diffusivity -- and expressing the temperatures relative to initial temperature $T_0$, i.e.
\begin{equation}
\widetilde{T}_{out}=T_{out}-T_0,\quad {\widetilde{T}}_{in}= T_{in} - T_0,\quad {\widetilde{T}}_g = T_g - T_0,
\end{equation}
translates \eqref{eqsimplified} into an ODE for the rescaled output $\widetilde{T}_{out}$ of the LTI system \eqref{semi-discrete}, i.e.
\begin{equation}
\frac{d\widetilde{T}_{out}}{dt} = A\widetilde{T}_{out} + \widetilde{T}_{in}b_1 + \widetilde{T}_{g}b_2,\quad A = -\frac{\dot{m}C_{p_w}+h_w\mathcal{A}}{m C_{p_w}},\quad b_1 = \frac{\dot{m}}{m},\quad b_2 = \frac{h_w\mathcal{A}}{m C_{p_w}},
\label{eqcorrelation-simplified0}
\end{equation}
with initial condition $\widetilde{T}_{out}(0) = {\widetilde{T}}_0=0$. This relation has an analytical solution given by
\begin{equation}
\widetilde{T}_{out}(t) = U(t)\widetilde{T}_\infty = \widetilde{T}_{in}F_1(t) + \widetilde{T}_g F_2(t),
\label{eqcorrelation-simplified}
\end{equation}
where $U(t) = 1 - \exp(At)$ is the evolution operator that evolves $\widetilde{T}_{out}$ from initial state $\widetilde{T}_{out}(0) = 0$ to
asymptotic state $\widetilde{T}_\infty = -(\widetilde{T}_{in}b_1 + \widetilde{T}_{g}b_2)/A$. (Note that $A<0$ ensures convergence on $\widetilde{T}_\infty$ with a characteristic time $\tau = -1/A$.) The second equation
in \eqref{eqcorrelation-simplified} expresses the rescaled output $\widetilde{T}_{out}$ in terms of the rescaled inputs $\widetilde{T}_{in}$ and $\widetilde{T}_{g}$ and defines the
sought-after input-output relation for the simplified system in \Cref{pipe-simplified}, with
\begin{equation}
F_1(t) = -\frac{b_1 U(t)}{A} = \frac{\dot{m}C_{p_w}\left(1-\exp(At)\right)}{\dot{m}C_{p_w}+h_w\mathcal{A}},\quad
F_2(t) = -\frac{b_2 U(t)}{A} = \frac{h_w\mathcal{A}\left(1-\exp(At)\right)}{\dot{m}C_{p_w}+h_w\mathcal{A}},
\label{eqcorrelation-simplified2}
\end{equation}
the transfer functions that describe the response of the system to the input \cite{ControlTheoryEng}.

For the DH pipe segment shown in \Cref{sketch-pipe} the situation is more complicated compared to the above simplified problem yet in essence the same input-output relation as \eqref{eqcorrelation-simplified} follows from the LTI form \eqref{semi-discrete} of the FOM. The latter namely has an analytical solution with a similar structure as \eqref{eqcorrelation-simplified}, i.e.
\begin{equation}
\widetilde{\xvec{T}}(t) = \xvec{U}(t)\widetilde{\xvec{T}}_\infty = \widetilde{T}_{in}\xvec{f}_1(t) + \widetilde{T}_g \xvec{f}_2(t),
\label{Analytical1}
\end{equation}
with tildes again indicating temperatures relative to the uniform initial temperature $T_0$ according to \eqref{BCs}. Here the evolution
operator and asymptotic state are given by matrix $\xvec{U}(t) = \xvec{I}-e^{(\xvec{A}t)}$ and vector $\widetilde{\xvec{T}}_\infty = -\xvec{A}^{-1}(\widetilde{T}_{in}\xvec{b}_1 + \widetilde{T}_{g}\xvec{b}_2)$, respectively, and vectors
\begin{equation}
\xvec{f}_1(t) = -\xvec{U}(t)\xvec{A}^{-1}\xvec{b}_1,\quad\quad \xvec{f}_2(t) = -\xvec{U}(t)\xvec{A}^{-1}\xvec{b}_2,
\label{Analytical2}
\end{equation}
define the corresponding transfer functions between the inputs $(\widetilde{T}_{in},\widetilde{T}_{g})$ and temperature field $\widetilde{\xvec{T}}$ (refer to \ref{appsolution}
for a detailed derivation).
The outlet temperature follows from
%
\begin{equation}
\widetilde{T}_{out}(t) = \xvec{C} \widetilde{\xvec{T}}(t),
\label{Output}
\end{equation}
with $\xvec{C}$ the standard output matrix for dynamical systems \cite{ControlTheoryEng}, and readily yields
\begin{equation}
\widetilde{T}_{out}(t) = \widetilde{T}_{in}F_1(t) + \widetilde{T}_g F_2(t),\quad
F_1(t) = \xvec{C} \xvec{f}_1(t),\quad
F_2(t) = \xvec{C} \xvec{f}_2(t),
\label{eqrelation}
\end{equation}
with transfer functions $F_{1,2}(t)$, as input-output relation for the pipe segment in \Cref{sketch-pipe}.
This is indeed the same as input-output relation \eqref{eqcorrelation-simplified} for the simplified problem. The difference is mainly technical in that the transfer functions
are for the latter available in the explicit form \eqref{eqcorrelation-simplified2}.

Relation \eqref{eqrelation} forms the backbone for the ROM to be developed in the remainder of \Cref{ROM} and, by deriving from the generic semi-discrete model \eqref{semi-discrete}, holds (for reasons given before) for a wide range of systems beyond the 1D pipe configuration of \Cref{numericalmethods}. The only additional constraint for validity of the resulting ROM is a uniform initial condition $T_0$ (\ref{appsolution}).

\subsection{Construction of temperature evolution from sequences of unit step responses}\label{UnitResponse}

The leading and trailing terms on the RHS of input-output relation \eqref{eqrelation} describe the response of the output $\widetilde{T}_{out}$ to step-wise changes in the inlet temperature by amount $\widetilde{T}_{in}$ and the ground temperature by amount $\widetilde{T}_g$, respectively, at time $t=0$. Linearity of the system admits generalisation of this single-step response to an arbitrary sequence of step-wise changes at arbitrary time levels $t_j$ \cite{ControlTheoryEng}. Introduce to this end the Heaviside function $H\left(t\right)$ \citep{Kreyszig2008}, i.e.
\begin{equation}
H\left(t\right) = \left\{
\begin{aligned}
0 && {t \leq 0}\\
1 && {t > 0}\\
\end{aligned}\right.,
\label{Heaviside}
\end{equation}
which allows expression of a sequence of $K_1$ step-wise changes of the inlet temperature from ${\widetilde{T}}_{in,j-1}$ to ${\widetilde{T}}_{in,j}$ at time levels $t_j$, with $t_0=0$, in the functional form
\begin{eqnarray}\label{eqTin}
{\widetilde{T}}_{in}\left(t\right) = H(t)\widetilde{T}_{in,0} - H(t-t_1)\widetilde{T}_{in,0} + H(t-t_1)\widetilde{T}_{in,1} + \cdots = \sum_{j=0}^{K_1}H\left(t - t_j\right)\Delta\widetilde{T}_{in,j},
\end{eqnarray}
with $\Delta\widetilde{T}_{in,0} = \widetilde{T}_{in,0}$ and $\Delta\widetilde{T}_{in,j} = \widetilde{T}_{in,j} - \widetilde{T}_{in,j-1}$ for $j\geq 1$, and similarly for $K_2$ step-wise changes of the ground temperature from ${\widetilde{T}}_{g,j-1}$ to ${\widetilde{T}}_{g,j}$ at (different) time levels $t_j$ through
\begin{equation}\label{eqTg}
{\widetilde{T}}_{g}\left(t\right) =  \sum_{j=0}^{K_2}{H\left(t - t_j\right)\Delta\widetilde{T}}_{g,j},
\end{equation}
with $\Delta\widetilde{T}_{g,0} = \widetilde{T}_{g,0}$ and $\Delta\widetilde{T}_{g,j} = \widetilde{T}_{g,j} - \widetilde{T}_{g,j-1}$ for $j\geq 1$.
The step-wise changes in inlet and ground temperatures according to \eqref{eqTin} and \eqref{eqTg} yield via the beforementioned linearity of
the input-output relation a temperature response at the pipe outlet following
\begin{equation}
\widetilde{T}_{out}(t) =  \sum_{j=0}^{K_1}G_1\left(t - t_j\right)\Delta \widetilde{T}_{in,j} + \sum_{j=0}^{K_2}G_2\left(t - t_j\right)\Delta\widetilde{T}_{g,j},
\label{eqTout}
\end{equation}
where functions
\begin{equation}
G_{1,2}(t - t_j) = H(t - t_j)F_{1,2}(t - t_j),
\label{eqTout2}
\end{equation}
effectively activate the transfer functions $F_{1,2}(t)$ at each time level $t_j$.

Input-output relation \Cref{eqTout} can to good approximation also describe the response to continuous changes in inlet and ground temperatures upon using sufficiently short and constant time intervals $\Delta t_j = t_j - t_{j-1} = \Delta t$. Consider to this end the formal solution for output \eqref{Output} in case of a generic unsteady input $\xvec{g}(t)$, which is given by the convolutions
\begin{eqnarray}
\widetilde{T}_{out}(t) = \xvec{C}\widetilde{\xvec{T}}(t) = \int_0^t \xvec{Q}(t-\xi)\xvec{g}(\xi)d\xi = \int_0^t q_1(t-\xi)\widetilde{T}_{in}(\xi)d\xi + \int_0^t q_2(t-\xi)\widetilde{T}_{g}(\xi)d\xi,
\label{eqTout3}
\end{eqnarray}
with $\xvec{Q}(t) = [q_1\;q_2] = \xvec{C}\exp(\xvec{A}t)\xvec{B}$ and $q_{1,2} = \xvec{C}\exp(\xvec{A}t)\xvec{b}_{1,2}$ \cite{ControlTheoryEng}.
Integration by parts yields
\begin{eqnarray}
\widetilde{T}_{out}(t) &=&
\int_0^t F_1(t-\xi) \widetilde{T}_{in}'(\xi) d\xi + \int_0^t F_2(t-\xi) \widetilde{T}_{g}'(\xi) d\xi\nonumber\\
&=& \int_0^t G_1(t-\xi) \widetilde{T}_{in}'(\xi) d\xi + \int_0^t G_2(t-\xi) \widetilde{T}_{g}'(\xi) d\xi,
\label{eqTout4}
\end{eqnarray}
where $\widetilde{T}_{in,g}' = d\widetilde{T}_{in,g}/d t$ and using $F_{1,2}(t-\xi)=G_{1,2}(t-\xi)$ due to $t-\xi\geq 0$, and coincides with \eqref{eqTout} in the limit of vanishing $\Delta t$. Here $q_{1,2}(t)$ relate to the beforementioned transfer functions via
$F_{1,2}(t) = \int_0^t q_{1,2}(t-\xi)d\xi$. This demonstrates that \eqref{eqTout} indeed admits approximation of the temperature response at the pipe outlet to arbitrary inputs for sufficiently small $\Delta t$, which is essential for modelling networks of pipe segments.

\subsection{Identification  and compact representation of transfer functions}\label{FunctionIdentify}

Transfer functions $F_{1,2}$ via \eqref{eqTout2} form the backbone of the input-output relation \eqref{eqTout} and are formally defined by the system matrices of the
LTI system \eqref{semi-discrete} according to \eqref{Analytical2} and \eqref{eqrelation}. However, rather than construction them in this cumbersome way, a more practical alternative exists in identifying $F_{1,2}$ from data.
This relies on two special cases for the temperature evolution:
\begin{itemize}
\item \textit{Case 1}: \Cref{eqcorrelation-simplified} becomes $\widetilde{T}_{out}\left(t\right) = F_1(t)$ when $\widetilde{T}_g = 0$ and $\widetilde{T}_{in} = H(t)$.
\item \textit{Case 2}: \Cref{eqcorrelation-simplified} becomes $\widetilde{T}_{out}\left(t\right) = F_2(t)$ when $\widetilde{T}_g = H(t)$ and $\widetilde{T}_{in} = 0$.
\end{itemize}
The sought-after transfer functions $F_1(t)$ and $F_2(t)$ can thus directly be identified from the two responses as described in {\it Case 1} and {\it Case 2}. They define the so-called ``unit step responses'' (i.e. they correspond with changes in the inputs from zero to unity at $t=0$ \cite{ControlTheoryEng}) and can be determined from both computational and experimental data. Here computational data from simulations by the FOM of \Cref{FOM} will be employed by the procedure elaborated below. However, it must be stressed that essentially the same approach can be adopted for data obtained by any computational method, laboratory
experiments or even field measurements for any pipe configuration governed by laws of physics that admit expression in the generic LTI form \eqref{semi-discrete}.

\begin{figure}[h!]
\centering
\includegraphics[width=7cm, height=4.5cm]{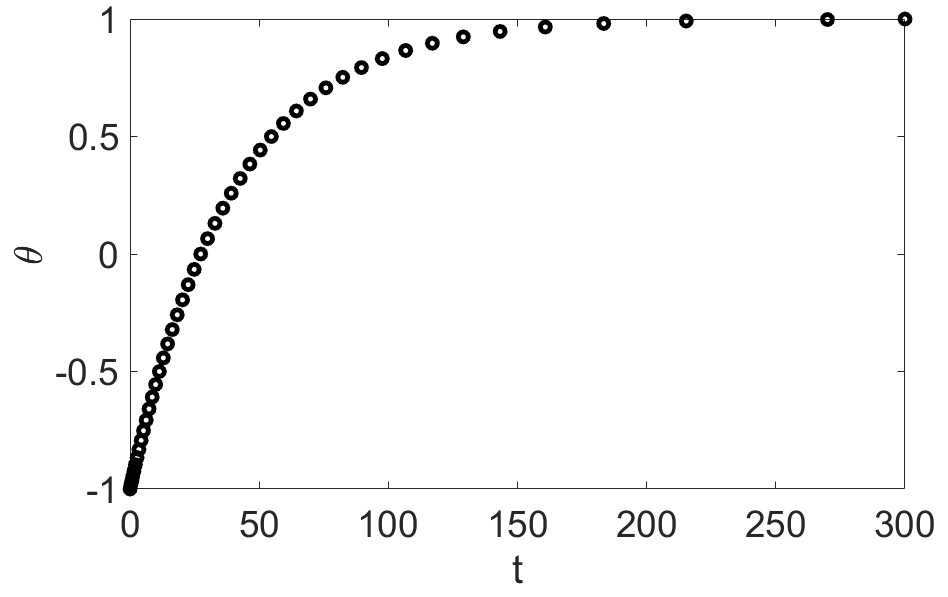}
\caption{Mapping $\theta(t)$ following \eqref{eqCheby} between spectral space $-1\leq \theta \leq 1$ and physical time $t\geq 0$ shown
in interval $0\leq t \leq 300$ for $t_{max}=300$ and $N=48$.}
\label{tvstheta}
\end{figure}

The identification procedure hinges on semi-analytical expression of the transfer functions in a basis of Chebyshev polynomials $\phi_n$ \citep{Canuto2012}:
\begin{equation}\label{eqCheby}
F(t) = \sum_{n=0}^{N}\hat{F}_n\phi_n\left(\theta(t)\right),\quad\quad \phi_n(\theta) = \cos(n\arccos\theta),\quad\quad\theta(t) = 1-2e^{-t/\tau},
\end{equation}
where $\phi_n(\theta)$ is the $n$-th order Chebyshev polynomial defined on the spectral space $\theta\in[-1,1]$ and $\hat{F}_n$ is the Chebyshev spectrum. Transformation $\theta(t)$ maps the spectral space
onto the physical semi-infinite time interval $t\in \left[0,\infty\right)$ and thus ensures that expansion \eqref{eqCheby} can capture the entire step response from initial state $F_{1,2}(0)$ to asymptotic state $F_{1,2}^\infty = \lim_{t\rightarrow\infty}F_{1,2}(t) = F_{1,2}|_{\theta=1}$ \citep{Canuto2012}.
Time constant $\tau$ is determined via relation
%
$\tau = -t_{ref}/\log\left[(1-\theta_{ref})/2\right]$
%
using a finite reference time $t_{ref}$ sufficiently far into the transient towards $F_{1,2}^\infty$ and $\theta_{ref}<\max\theta = 1$ as corresponding reference in spectral space.
Here $\theta_{ref} = \cos(\pi/N)$ for reasons elucidated below and $t_{ref} = \alpha t_{max}$, with $\alpha = 0.9$ and $t_{max}$ the final temperature of the data set. \Cref{tvstheta} gives the
resulting transformation $\theta(t)$ for $t_{max}=300$ and $N=48$, yielding $\tau = 39.5$.

Representing the transfer functions in orthogonal polynomials, here Chebyshev polynomials according to \Cref{eqCheby}, has important advantages. First, these expansions exhibit so-called ``spectral convergence''
and thus enable accurate description of $F_1(t)$ and $F_2(t)$ at any time $t$ by Chebyshev expansion \eqref{eqCheby} with an order $N$ that is substantially lower than a standard expansion based on e.g. piece-wise
linear interpolation \citep{Canuto2012}. The spectral approach thus effectively yields a reduced-order model of the input-output relation compared to a conventional representation and thereby contributes to
a further reduction of the model size. Second, for given Chebyshev spectrum ${\hat{F}}_n$, expansions \eqref{eqCheby} constitute global semi-analytical expressions in terms of basic trigonometric functions and thus admit evaluation of function values at
arbitrary time levels $t$ as well as easy incorporation in larger system models and exact performance of mathematical operations such as differentiation and integration for e.g. control purposes. This approach thereby provides a compact and mathematically sound way to describe the input-output relation in \Cref{eqTout}.

The Chebyshev spectrum ${\hat{F}}_n$ in \eqref{eqCheby} is determined via the discrete Chebyshev transform
\begin{equation}\label{ChebySpectral}
{\hat{F}}_n=\frac{1}{\gamma_n}\sum_{k=0}^{N}F(\theta_k)\phi_n(\theta_k)w_k, \;\;\;\;\; \gamma_n = \left\{
\begin{aligned}
\pi && n=0,N\\
\frac{\pi}{2}  && 1 \leq n \leq N-1\\
\end{aligned} \right.,
\end{equation}
with $\theta_k \in[-1,1]$ and $w_k$ the quadrature points and weights, respectively, given by
\begin{equation}
\theta_k=cos\frac{\pi k}{N},\;\;\;\;\; w_k = \left\{
\begin{aligned}
\frac{\pi}{2N} && k=0,N\\
\frac{\pi}{N}  && 1 \leq k \leq N-1\\
\end{aligned} \right.
\end{equation}
according to the Chebyshev-Gauss-Lobatto procedure \citep{Canuto2012}. The function values $F\left(\theta_k\right)$ required for Chebyshev
transform \eqref{ChebySpectral} constitute the training data for expansion \eqref{eqCheby} and are determined from interpolation of FOM training data of {\it Case 1} and {\it Case 2} defined before on the discrete time levels $t_k = -\tau\log\left[(1-\theta_k)/2\right]$ corresponding with transformation $\theta(t)$ in \eqref{eqCheby}. Note that reference $\theta_{ref} = \cos(\pi/N)$ for the evaluation of time constant $\tau$ coincides with the first internal point $\theta_1$ relative to interval boundary $\theta_0 = 1$.

\section{Performance and validation of the ROM}\label{discussion}

\subsection{Experimental validation of the FOM}

The performance of the FOM described in \Cref{FOM} is validated by using experimental data from \citep{Sartor2017}. To this end they realized an experimental set-up consisting of a steel pipe of length $L=39 m$,
an inner diameter of $D=52.48 mm$ and a thickness of $3.91 mm$, yielding $r_w = 26.24 mm$ and $r_s = 30.15 mm$ for the corresponding radii in \Cref{sketch-pipe}. The pipe is insulated with Tubolit 60/13 (thickness $13 mm$) to mimic the insulation layer of DH pipes and has a radius of $r_i = 43.15 mm$. The casing layer shown in \Cref{sketch-pipe} is not included in the set-up, however, and the corresponding energy balance \Cref{eqcasing} is therefore eliminated from the numerical model and temperature $T_c$ in the trailing term in \eqref{eqinsu} is substituted by the ambient temperature $T_g$. Another essential difference between real DH pipes and the one used in the experimental set-up is that DH pipes are buried underground while the insulated pipe in the set-up is exposed to ambient air. However, these differences between actual DH pipes and the experimental configuration are non-essential, since the overall system behavior remains qualitatively the same, i.e. convective heat exchange between water and inner pipe wall, conduction through multiple solid layers and convective heat exchange between the outer layer and the environment. Hence the experimental set-up is well-suited for validation of the FOM (refer to \citep{Sartor2017} for further details on this set-up and the corresponding data).
\Cref{TableExpValidation} gives the relevant model parameters for the experimental validation of the FOM.
\begin{table}[h!]
\begin{centering}
\caption{Model parameters for experimental validation of the FOM.}
\label{TableExpValidation}
\begin{tabular}{llll}\\\hline
Element & \begin{tabular}[c]{@{}l@{}}$\rho$ ($kg/m^3$)\end{tabular} & \begin{tabular}[c]{@{}l@{}}$C_p (J/(kg K)$\end{tabular} & \begin{tabular}[c]{@{}l@{}}$k (W/(m K)$\end{tabular} \\ \hline
Steel pipe & 7800 & 480 & 45 \\
Tubolit 60/13 & 25 & 2450.7 & 0.04\\\hline
\end{tabular}\\
\end{centering}
\end{table}

\begin{figure}[htbp]
\begin{subfigure}[b]{0.45\textwidth}
\centering
\includegraphics[width=5.5cm,height=4.5cm]{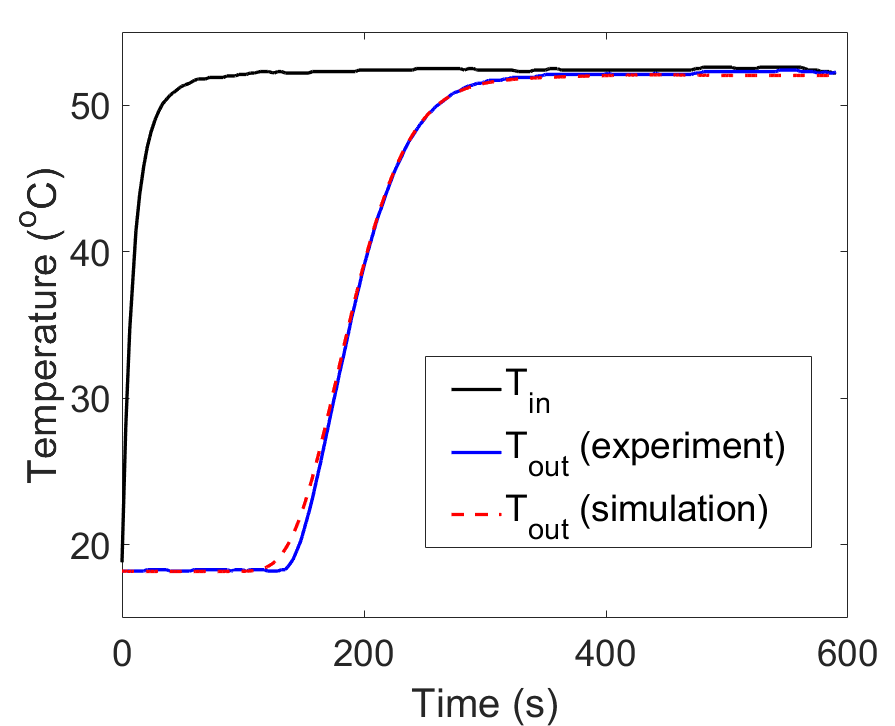}
	\caption{$u=0.27 m/s$}
	\label{rig-2}
\end{subfigure}
\begin{subfigure}[b]{0.45\textwidth}
	\centering
\includegraphics[width=5.5cm,height=4.5cm]{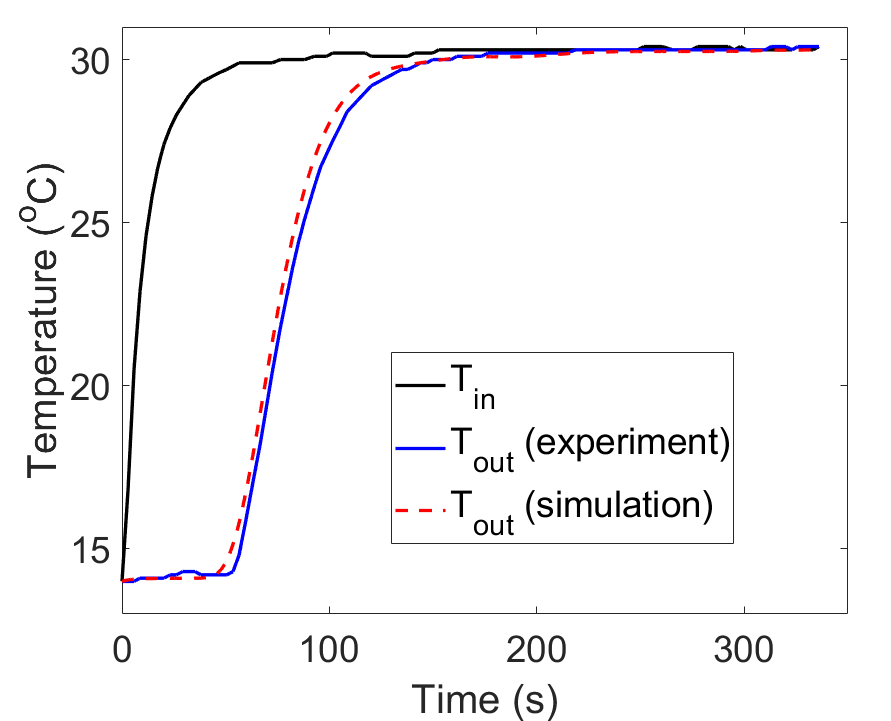}
	\caption{$u=0.74 m/s$}
	\label{rig-3}
\end{subfigure}
\hfill
\begin{subfigure}[b]{0.45\textwidth}
	\centering
\includegraphics[width=5.5cm,height=4.5cm]{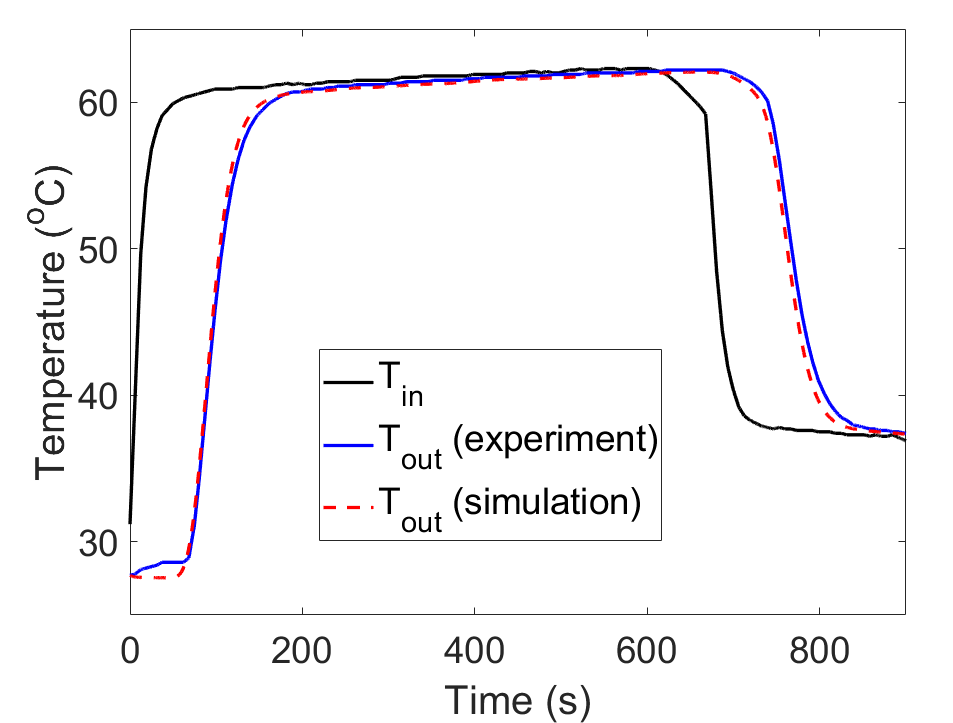}
	\caption{$u=0.57 m/s$}
	\label{rig-5}
\end{subfigure}
\hfill
\begin{subfigure}[b]{0.45\textwidth}
	\centering
\includegraphics[width=5.5cm,height=4.5cm]{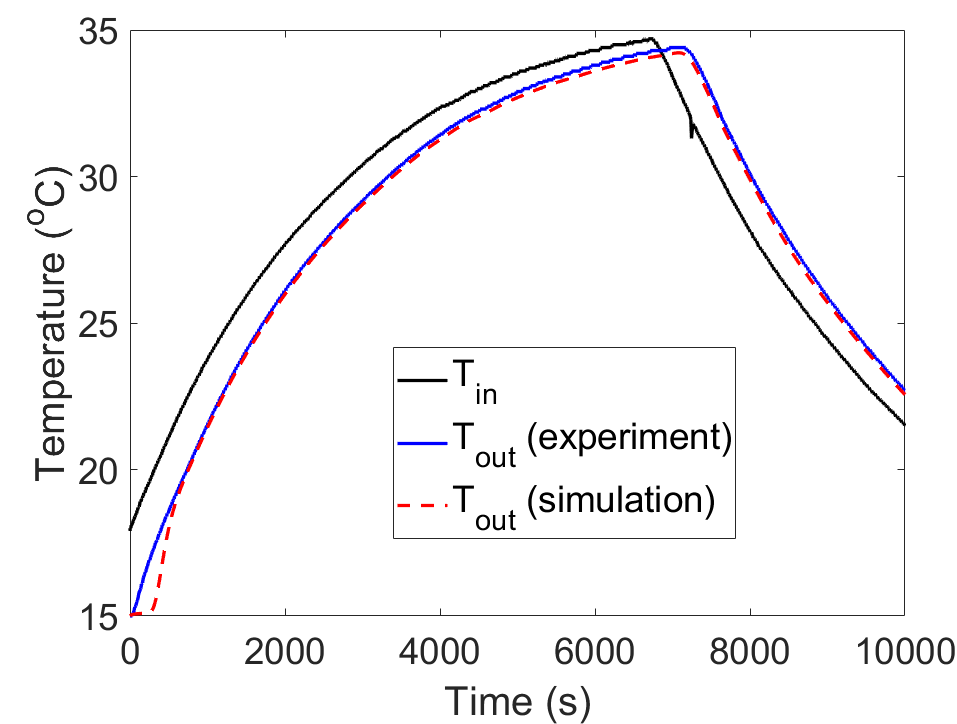}
	\caption{$u=0.11 m/s$}
	\label{rig-6}
\end{subfigure}

\caption{Evolution of the outlet temperature $T_{out}$ for given inlet temperature $T_{in}$ according to laboratory experiments from \citep{Sartor2017} versus FOM simulations for flow velocities $u$ as indicated.}

\label{FOMvsTest}
\end{figure}

A mesh-dependence analysis shows that a grid size of $\Delta x = 0.6 m$ is required for accurate results; the corresponding time step $\Delta t = \Delta t_{max}$ is determined according to stability
criterion \eqref{stability}. Comparison of the outlet temperature $T_{out}$ of the experiments versus simulations for given inlet profile $T_{in}$ at flow velocities ranging from
$u=0.11 m/s$ to $u=0.74 m/s$ is shown in \Cref{FOMvsTest}. This reveals a close agreement between measured and predicted response $T_{out}$ to both the step-like changes of $T_{in}$ in \Cref{rig-2}--\ref{rig-5} and
the gradual change of $T_{in}$  in \Cref{rig-6}. This demonstrates that the FOM admits sufficiently accurate prediction of the temperature distribution in DH networks to provide reliable training data for the ROM.

\subsection{Identification of the transfer functions from training data}\label{identifications}

The transfer functions $F_1(t)$ and $F_2(t)$ are identified from the unit step responses according to \Cref{FunctionIdentify} as simulated
with the FOM for a DH pipe segment following \Cref{sketch-pipe} with length $L=100m$, nominal diameter $D=25mm$ (DN25) and velocity $u=1.5m/s$. Further dimensions and properties of the insulated pipe are in accordance with the Logstor standard\footnote{Product catalogue for district energy (version 2020.03): \url{https://www.logstor.com/media/6506/product-catalogue-uk-202003.pdf}} and yield the relevant model parameters
for the identification of the transfer functions from FOM data shown in \Cref{TableSystemID}. These parameter values are used for all the test cases in \Cref{identifications} and \Cref{section3.3} unless stated otherwise.
\begin{table}[h!]
\centering
\caption{Model parameters for identification of the transfer functions from FOM data.}
\label{TableSystemID}
\begin{tabular}{llllll}\\\hline
Element & \begin{tabular}[c]{@{}l@{}}$r (mm)$\end{tabular} & \begin{tabular}[c]{@{}l@{}}Thickness\\ ($mm$)\end{tabular} & \begin{tabular}[c]{@{}l@{}} $\rho$ ($kg/m^3$)\end{tabular} & \begin{tabular}[c]{@{}l@{}} $C_p (J/(kg K)$\end{tabular} & \begin{tabular}[c]{@{}l@{}}$k (W/(m K)$\end{tabular} \\\hline
Water & 14.25 & N/R & 996.7 & 4066.7 & 0.605 \\
Steel pipe & 16.85 & 2.6 & 7900 & 502.5 & 51 \\
Insulation & 42 & 25.15 & 30 & 1400 & 0.027 \\
Casing & 45 & 3 & 944 & 2250 & 0.43 \\
Soil & N/R & N/R & N/R & N/R & 1.6\\\hline
\end{tabular}
\end{table}

\begin{figure}[htbp]
\centering
	\begin{subfigure}[b]{0.48\textwidth}
\includegraphics[width=6cm,height=4.5cm]{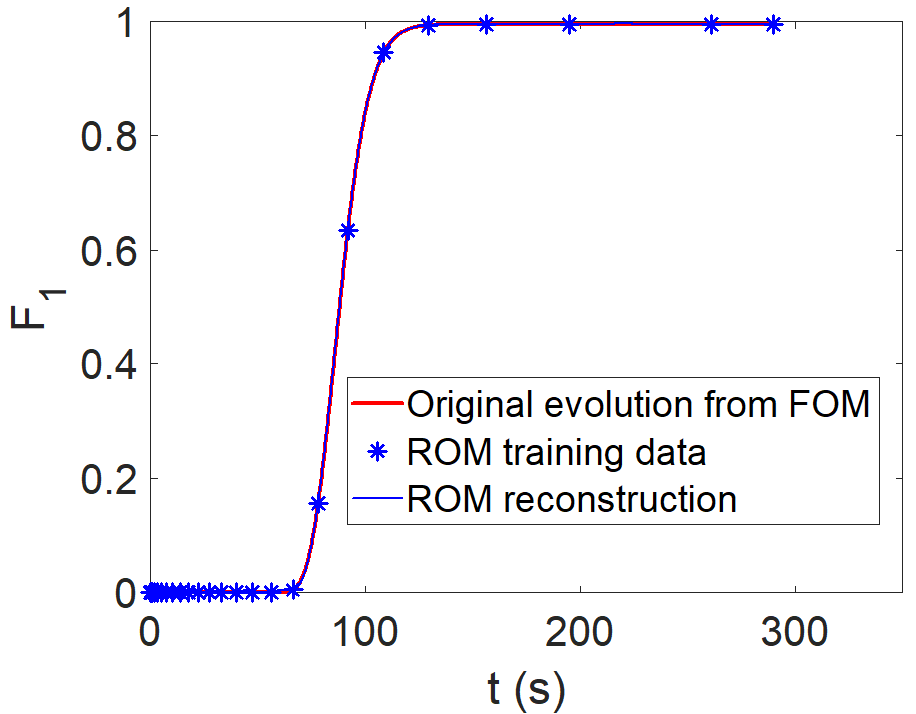}
	\caption{}\label{unit1}
	\end{subfigure}
	\hfill
	\begin{subfigure}[b]{0.48\textwidth}
\includegraphics[width=6cm,height=4.8cm]{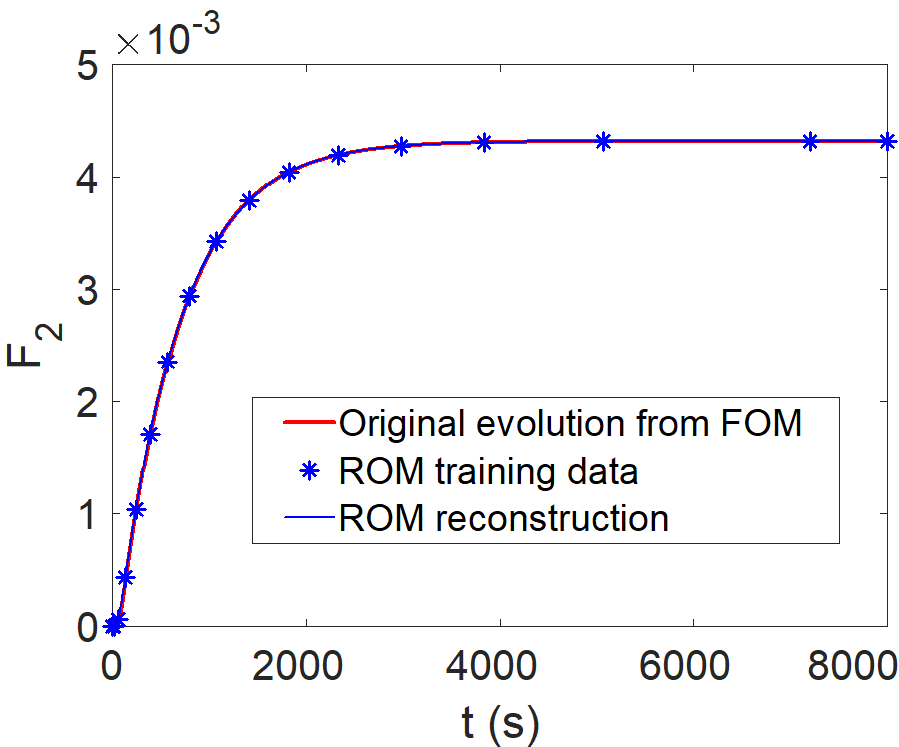}
	\caption{}\label{unit2}
	\end{subfigure}
\caption{Identification of transfer functions $F_1$ (panel a) and $F_2$ (panel b) of the ROM by the procedure following \Cref{UnitResponse} from training data generated by the FOM.}
\label{fgUnitresponse}
\end{figure}

\Cref{fgUnitresponse} gives the transfer functions $F_1$ and $F_2$ determined from FOM data (red curves) by the procedure following cases 1 and 2 in \Cref{UnitResponse}. The blue stars
indicate the discrete time levels $t_k$ for generation of the ROM training data for the evaluation of the Chebyshev spectrum \eqref{ChebySpectral} and the blue curves show the ROM
reconstruction of the transfer functions via \eqref{eqCheby}. Here expansion orders $N=24$ and $N=16$ have been taken for $F_1$ and $F_2$, respectively, and the close agreement between original
FOM evolution and ROM reconstruction signifies an accurate approximation by the latter using a very small number of Chebyshev polynomials.
Since the FOM prediction concerns the full system involving both the fluid and the surrounding layers (i.e. steel pipe, insulation, casing), all relevant physical phenomena
such as thermal losses and thermal inertia are included in the unit step responses and, inherently, in the identified transfer functions $F_{1,2}$. The evolutions of the transfer functions
show a number of notable differences, though. Transfer function $F_1$ (\Cref{unit1}) exhibits a delayed response at around 66s yet reaches its equilibrium $F_1^\infty$ quickly
at around 130s; transfer function $F_2$ (\Cref{unit2}) , on the other hand, responds immediately yet reaches its equilibrium $F_2^\infty$ only at around 5000s.
Moreover, the equilibria differ substantially in magnitude, i.e. $F_1^\infty\approx 1$ versus $F_2^\infty\approx 4.2 \times 10^{-3}$. Thus the change in ${\widetilde{T}}_{out}$ caused by
${\widetilde{T}}_{g}$ is much smaller and slower than the one caused by ${\widetilde{T}}_{in}$, implying system dynamics that are dominated by the inlet temperature.

The transfer functions in fact offer important insight into the system behavior. The differences in sensitivity to both inputs can namely be (qualitatively) understood from the transfer functions \eqref{eqcorrelation-simplified2} of the simplified model \eqref{eqsimplified}. These functions evolve monotonically from $F_{1,2}(0)=0$ to $\lim_{t\rightarrow\infty}F_{1,2}(t) = F_{1,2}^\infty$, with $F_{1}^\infty = \dot{m}C_{p_w}/(\dot{m}C_{p_w}+h_w\mathcal{A})<1$ and $F_{2}^\infty = h_w\mathcal{A}/(\dot{m}C_{p_w}+h_w\mathcal{A})<1$ the corresponding equilibria satisfying $F_{1}^\infty+F_{2}^\infty=1$. The system properties
according to \Cref{TableSystemID} imply $h_w\mathcal{A}\ll\dot{m}C_{p_w}$ and, in consequence, a clear dominance of the water influx at the inlet (i.e. $F_{1}^\infty\approx 1$) over the heat exchange with
the environment (i.e. $F_{2}^\infty\ll 1$) in the dynamics of the outlet temperature $T_{out}$. The dynamics of the FOM are more intricate, as e.g. reflected in the delayed response to $T_{in}$ and the different
equilibration times for $F_1$ and $F_2$, yet the general trend is similar: a clear dominance of the water influx for basically the same reasons (i.e. $F_1^\infty\gg F_2^\infty$).

\begin{figure}[h!]
\centering
	\begin{subfigure}[b]{0.45\textwidth}
\includegraphics[width=5.5cm,height=4.2cm]{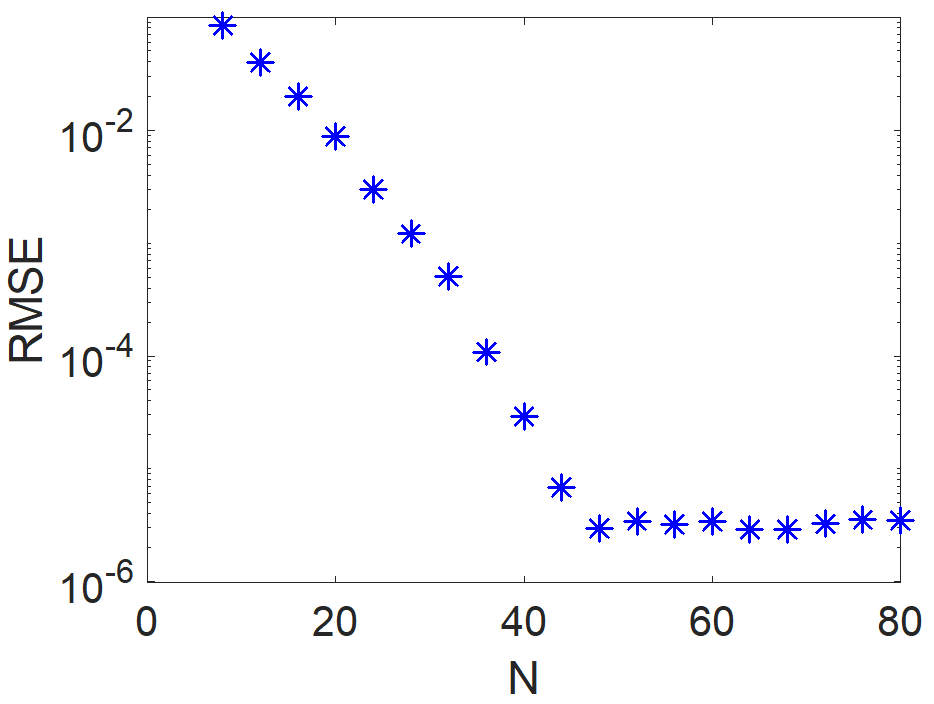}
	\caption{}\label{RMSE1}
	\end{subfigure}
	\hfill
	\begin{subfigure}[b]{0.45\textwidth}
\includegraphics[width=5.5cm,height=4.2cm]{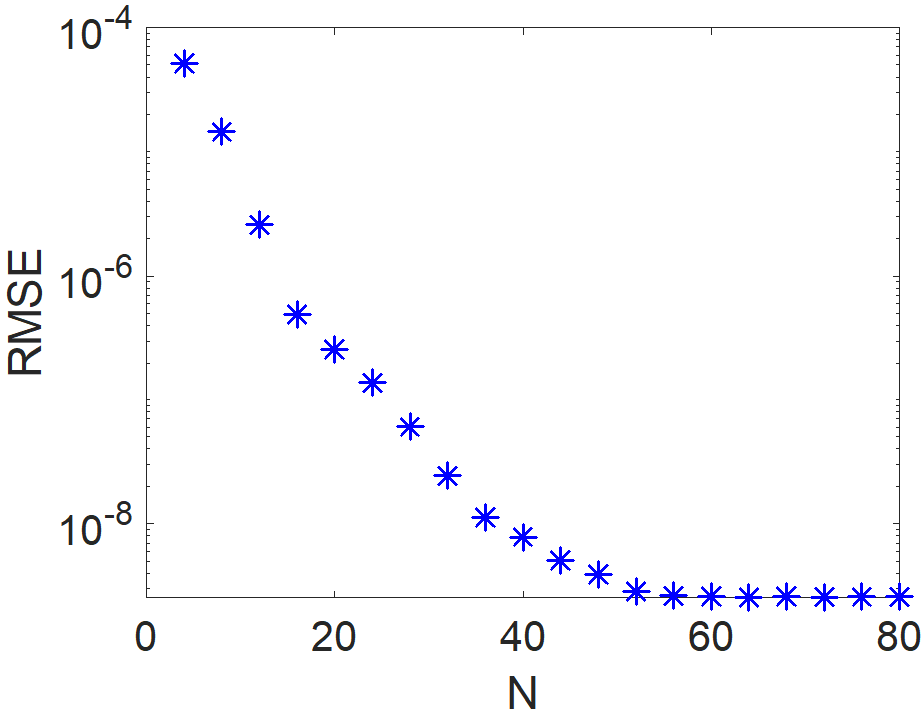}
	\caption{}\label{RMSE2}
	\end{subfigure}
\caption{RMSE between the ROM predictions $F_1$ (panel a) and $F_2$ (panel b) and the FOM benchmark as a function of the Chebyshev expansion order $N$.}
\label{RMSE}
\end{figure}

Visual inspection of \Cref{fgUnitresponse} suggests that the above expansion orders $N$ enable accurate capturing of the transfer functions $F_{1,2}$ by the ROM. A more systematic
determination of the appropriate $N$ may follow from the root mean square error (RMSE) between the ROM predictions $ROM_j$ of a given variable and the corresponding FOM benchmark $FOM_j$ at time levels $t_j$, i.e.
\begin{eqnarray}
RMSE = \sqrt{\frac{\sum_{j=1}^Q (FOM_j - ROM_j)^2}{K}}
\label{RMSEdef}
\end{eqnarray}
with $Q$ the total number of time steps in a simulation. (Here the solution of the ROM is projected onto the same time levels $t_j$ as in the FOM by linear interpolation.)
\Cref{RMSE} gives the $RMSE$ versus the expansion order $N$ for both transfer functions $F_{1,2}$ and reveals an exponential decrease with growing $N$ until saturation sets in at $RMSE\approx 2\times 10^{-6}$ for $F_1$ at $N=48$ and $RMSE\approx 2\times 10^{-9}$ for $F_2$ at $N=52$. (This difference of 3 orders of magnitude is consistent with the
difference in magnitude of $F_1^\infty$ and $F_2^\infty$.) Setting the desired level of accuracy for the ROM predictions at $RMSE = 2\times 10^{-6}$ yields $N=48$ and $N=16$
as adequate expansion orders for $F_1$ and $F_2$, respectively. These orders are used for the remainder of this study unless stated otherwise.

\begin{figure}[h!]
\centering
	\begin{subfigure}[b]{0.45\textwidth}
\includegraphics[width=5.8cm,height=4.5cm]{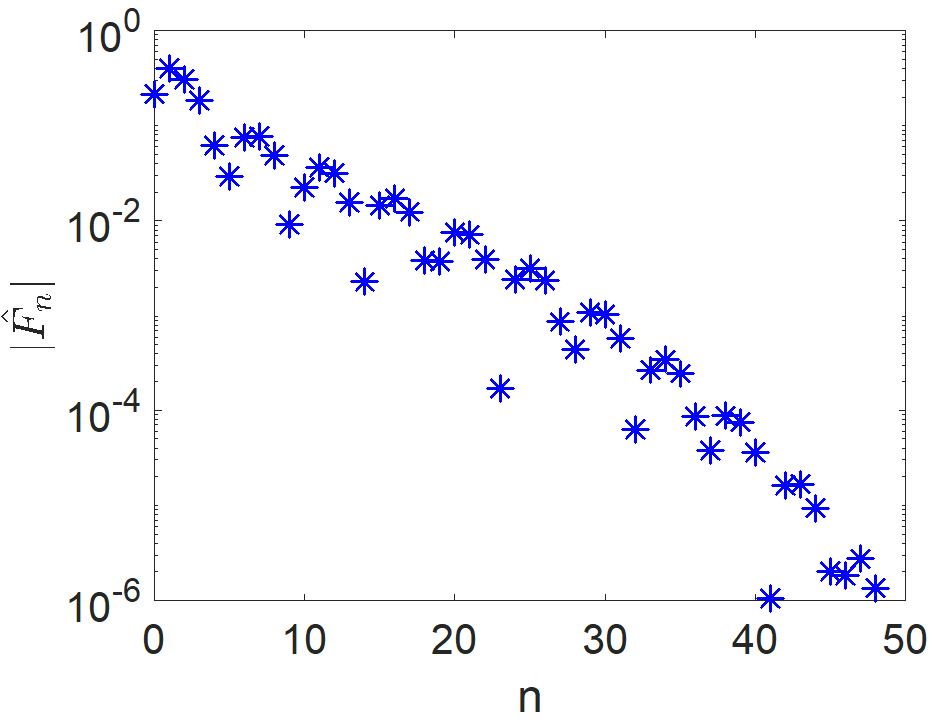}
	\caption{}\label{spectrum1}
	\end{subfigure}
	\hfill
	\begin{subfigure}[b]{0.45\textwidth}
\includegraphics[width=5.8cm,height=4.5cm]{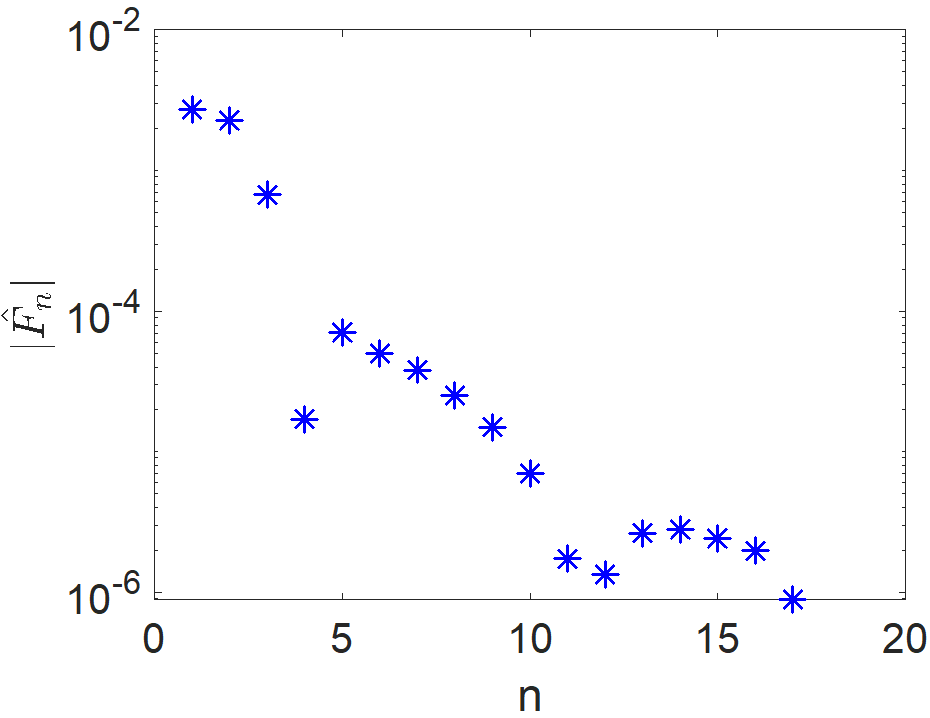}
	\caption{}\label{spectrum2}
	\end{subfigure}
\caption{Chebyshev spectrum ${\hat{F}}_n$ of the identified transfer functions $F_1$ (panel a) and $F_2$ (panel b) for expansion orders $N = 48$ and $N = 16$, respectively.}	
\label{spectra}
\end{figure}

The Chebyshev spectra ${\hat{F}}_n$ corresponding with the transfer functions $F_1$ and $F_2$ of \Cref{fgUnitresponse} for $N$ as determined above are shown in \Cref{spectra}. This clearly reveals the exponential decay of the magnitude of coefficients ${\hat{F}}_n$ with increasing polynomial order $n$ that characterizes spectral convergence of expansions in orthogonal polynomials for sufficiently large $N$ \cite{Canuto2012}. This spectral convergence is consistent with the exponential decay of the RMSE towards saturation in \Cref{RMSE} and constitutes further evidence of an accurate identification of the transfer functions from the FOM data with the chosen $N$. Important to note is that this degree of accuracy is attained for any time t, i.e. not just at the discrete time levels $t_k$ corresponding with the training data. This is a direct consequence of the convergence properties of the employed type of
polynomial expansions \citep{Canuto2012}. These properties admit reliable and efficient prediction of the outlet temperature $T_{out}$ in response to step-wise changes in $T_{in}$ and $T_g$ via the input-output relation \eqref{eqTout} at any time $t$.

\subsection{Response to arbitrary step-wise input profiles}\label{section3.3}

\begin{figure}[]
\centering
	\begin{subfigure}[b]{0.45\textwidth}
	\includegraphics[width=7cm,height=5cm]{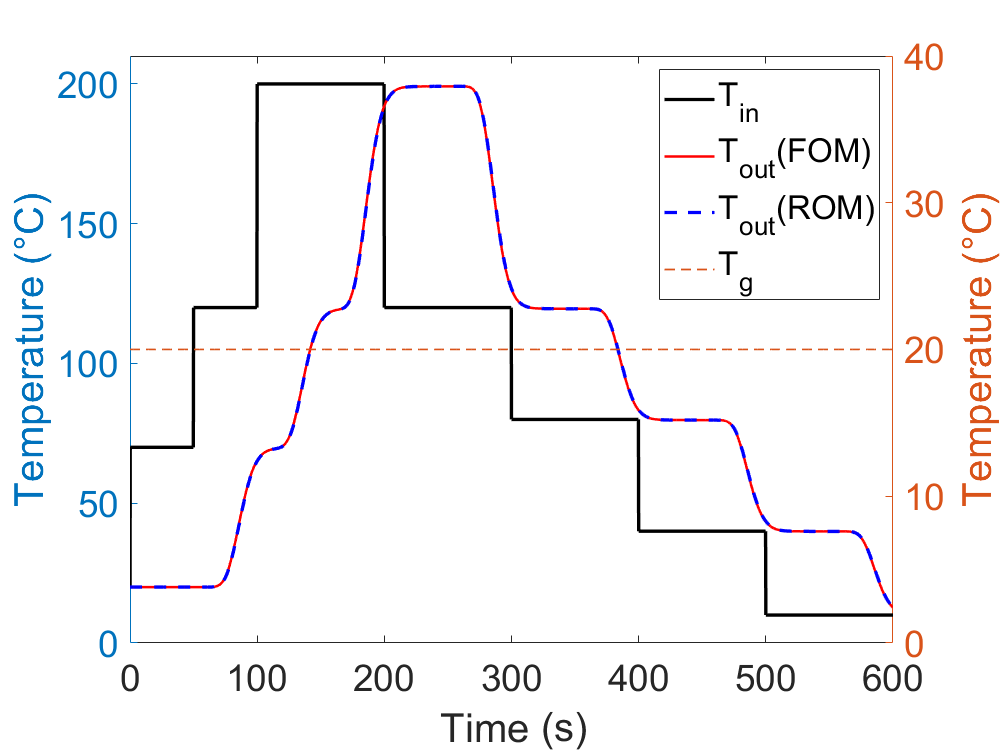}
	\caption{}\label{case1}
	\end{subfigure}
	\hfill
	\begin{subfigure}[b]{0.45\textwidth}
	\includegraphics[width=7cm,height=5cm]{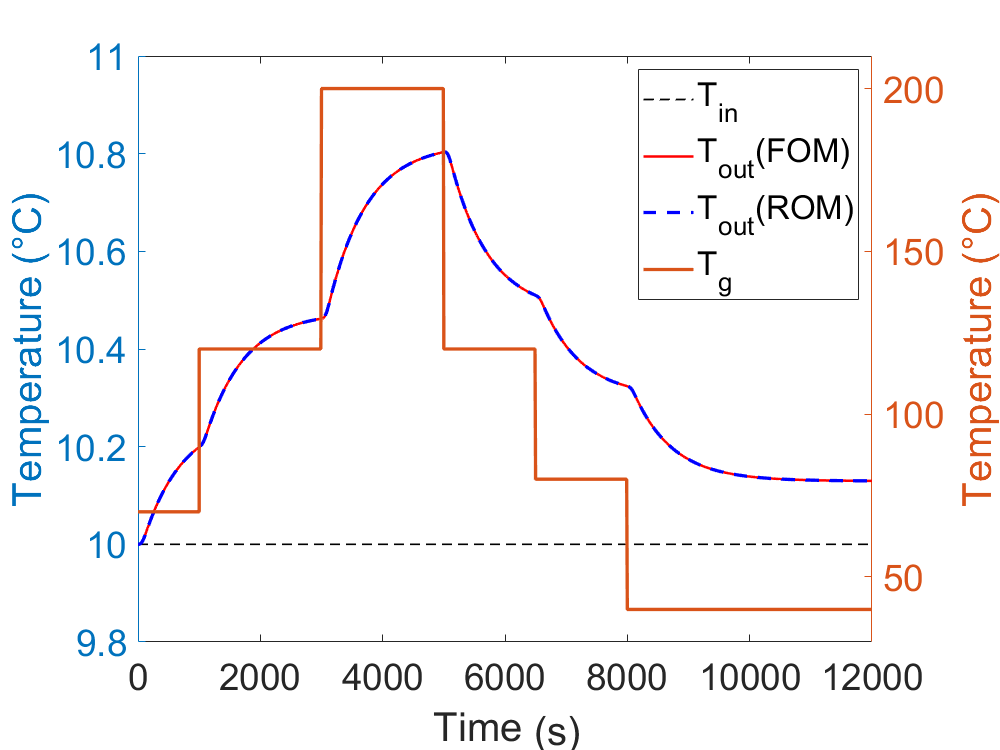}
	\caption{}\label{case2}
	\end{subfigure}

\caption{Response of $T_{out}$ to step-wise changes in either $T_{in}$ (panel a) or $T_{g}$ (panel b) according to ROM predictions versus FOM simulations.}

\label{cases}
\end{figure}

The performance and validity of the ROM are investigated further by determining the response of $T_{out}$ to multiple step-wise changes in the input. \Cref{cases} shows $T_{out}$ predicted by the
ROM (dashed blue; left axis) as a function of the given profiles for the inlet temperature $T_{in}$ (solid blue; left axis) and ground temperature $T_{g}$ (solid orange; right axis) versus $T_{out}$ according to the
FOM (solid red). Here only input $T_{in}$ (\Cref{case1}) or $T_{g}$ (\Cref{case2}) is varied while keeping the other input constant. The response of $T_{out}$ to changes in $T_g$ (\Cref{case2}) is
considerably slower and weaker compared to similar changes in $T_{in}$ (\Cref{case1}) due to the dramatically different sensitivity to former and latter input (and magnitude of the corresponding
transfer function $F_2$ and $F_1$, respectively) found in \Cref{identifications}. However, in both cases (again) a close agreement between
the ROM and FOM occurs and this further demonstrates that the response to step-wise changes in $T_{in,g}$ are accurately captured by transfer functions $F_{1,2}$.

Note that, for illustration and testing purposes, exaggerated profiles are used in \Cref{cases} for the variations in $T_{in}$ and $T_g$. In practice, variations in $T_g$ are much smaller and more gradual
due to the large thermal mass and inertia of the ground; variations in $T_{in}$ are usually also more moderate in new generations of DH networks due to lower heating demands of modern buildings and
better insulation of the piping.
Important furthermore is that, since the temperature responses for multiple step-wise changes are constructed by linear combination of the unit step responses according to input-output relation (\ref{eqTout}), thermo-physical
properties of the fluid must be assumed constant so as to ensure linearity of the thermal behavior and validity of this relation. For this purpose, constant fluid properties are evaluated using
the mean inlet temperature as a reference.

\begin{figure}[h!]
\centering
\includegraphics[width=8cm, height=5cm]{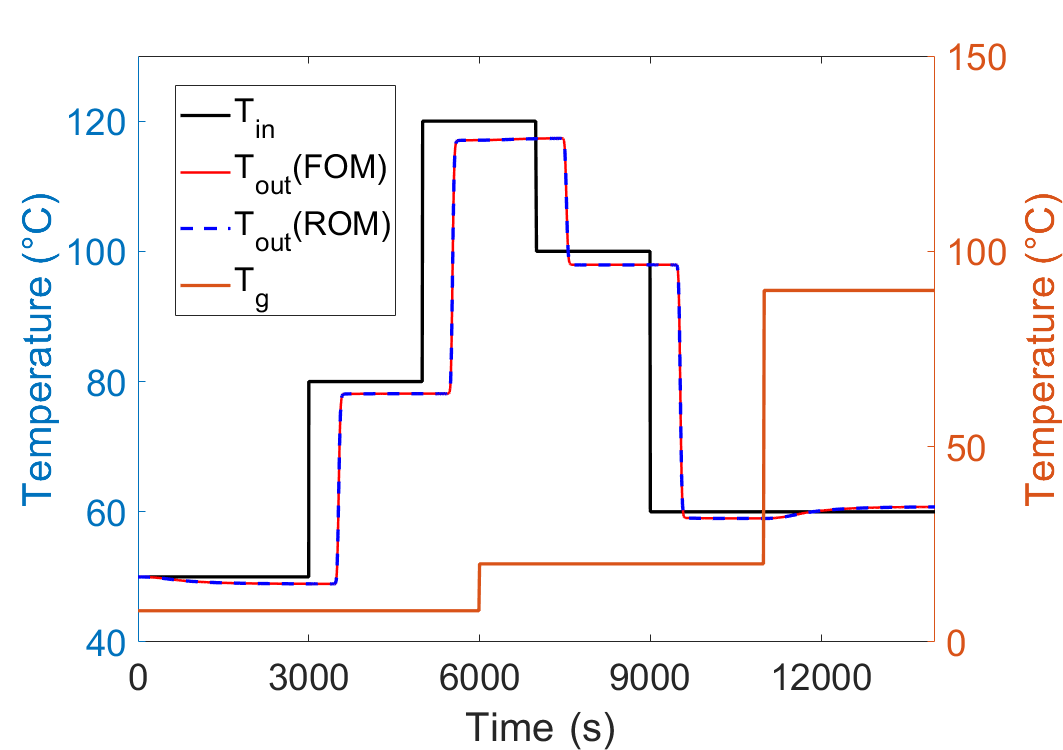}

\caption{Response of $T_{out}$ to step-wise changes in both $T_{in}$ and $T_{g}$ according to ROM predictions versus FOM simulations. Note: concerns pipe segment with length $L=600m$ instead of $L=100m$ in \Cref{cases}.}
\label{case3}
\end{figure}

\Cref{case3} shows the response of $T_{out}$ to combined step-wise changes in both $T_{in}$ and $T_g$ following the given profiles. Here a longer pipe of length $L=600m$ is used, since the impact of $T_g$ is
very weak compared to that of $T_{in}$ for the original pipe length (\Cref{case2}). (Other parameters remain as per \Cref{TableSystemID}.) From the simplified model \eqref{eqsimplified} it namely readily follows that
larger $L$ yields a larger area $\mathcal{A}$ and via an increased magnitude of transfer function $F_2$ in \eqref{eqcorrelation-simplified2} thus a stronger influence of $T_g$; similar boosting of the impact of
$T_g$ with larger pipe length occurs for the full system. (Maintaining an equal degree of accuracy necessitates an increase of the expansion orders for the identification of the transfer functions
to $N = 96$ for $F_1$ and $N = 24$ for $F_2$.)
Comparison of the evolutions according to ROM and FOM in \Cref{case3} again reveals a close agreement and thereby shows that the ROM can accurately predict the temperature response at the outlet of
a DH pipe to any sequence of step-wise changes in the inlet and ground temperatures. Note that again exaggerated input temperatures are considered for illustration and testing purposes.

The behavior in \Cref{case3} demonstrates that, despite considering a much longer pipe, a significantly different sensitivity of $T_{out}$ to both inputs remains. The first stage up to $t\approx 3000s$ concerns a relatively low inlet
temperature $T_{in} = 50^\circ C$ and ground temperature $T_g = 10^\circ C < T_{in}$ and results in a net thermal loss (i.e. a gradual decrease of $T_{out}$ over time). From $t=3000s$ onward
$T_{in}$ step-wise increases to $T_{in} = 120^\circ C$ and subsequently decreases via $T_{in} = 100^\circ C$ back
to $T_{in} = 60^\circ C$. This input clearly dominates the corresponding response of $T_{out}$ in that the latter closely shadows the profile of $T_{in}$; the additional step-wise increase of
$T_g$ to $30^\circ C$ at $t=8000s$, on the other hand, has no noticeable effect on $T_{out}$. Only the significant reduction of $T_{in}$ from $100^\circ C$ to $60^\circ C$ at $t=9000s$ in
combination with a substantial jump in $T_g$ from $30^\circ C$ to $90^\circ C$ at $t=11,000s$ yields a clear impact of $T_g$ on $T_{out}$ in that the latter gradually increases in time.

\section{Application of the ROM}\label{application}

The analysis below concerns practical application of the ROM for the operation of DH systems. This is illustrated and examined for two scenarios: fast simulation of the thermal behavior of a
small DH network for e.g. optimal scheduling (\Cref{section4.1}) and design of a closed-loop controller for user-defined temperature regulation of a DH system (\Cref{ControlConsumer}).

\subsection{Fast simulation of a small DH network}\label{section4.1}

\begin{figure}[h!]
\centering
	\begin{subfigure}[b]{0.45\textwidth}
\includegraphics[width=7cm,height=3.5cm]{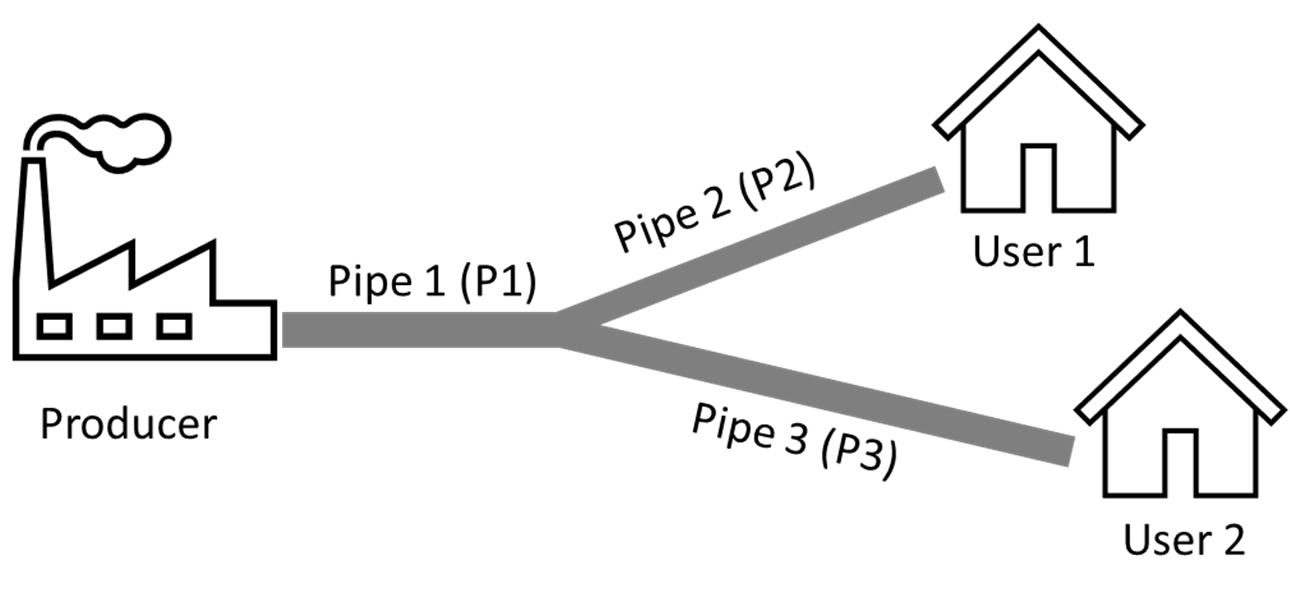}
	\caption{}
	\label{system1-1}
	\end{subfigure}
	\hfill
	\begin{subfigure}[b]{0.48\textwidth}
\includegraphics[width=7cm,height=3.5cm]{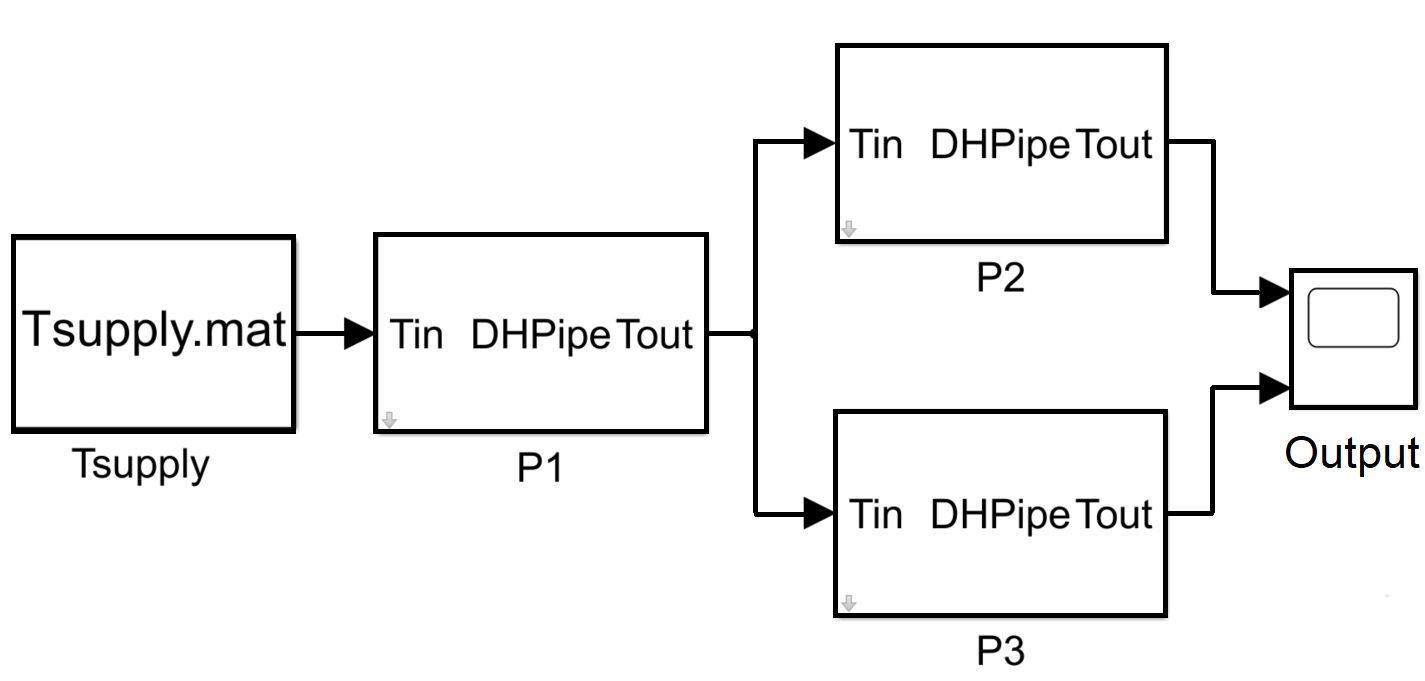}
	\caption{}
	\label{system1-2}
	\end{subfigure}
\caption{Small DH network \textit{System 1} consisting of one heat producer and two end users: (a) configuration; (b) {\tt Matlab-Simulink} block structure.}
\end{figure}

Considered is a small DH network with one energy producer and two consumers connected by pipes $P1$, $P2$ and $P3$ as sketched in \Cref{system1-1} and denoted \textit{System 1} hereafter. Water enters pipe $P1$ with a constant mass flux $\dot{m}_1 = 2.1884kg/s$ and a supply temperature $T_{supply}$ according to a typical daily variation for a DH system used in \citep{Jie2012} and is equally distributed over pipes $P2$ and $P3$
(i.e. $\dot{m}_2 = \dot{m}_3 =\dot{m}_1/2 = 1.0942 kg/s$). Pipe dimensions are specified in \Cref{tablesystem1} and yield flow velocities as indicated for given mass fluxes. The ground temperature is fixed at $T_g=10^\circ C$.\\\\

\begin{table}[h!]
\small
\centering
\caption{Model parameters of small DH network \textit{System 1}}
\label{tablesystem1}
\begin{tabular}{lllll}\hline
& $L(m)$ & $DN(mm)$ & $\dot m(kg/s)$ & $u(m/s)$ \\\hline
Pipe 1 (P1) & 200 & 40 & 2.1884 & 1.5  \\
Pipe 2 (P2) & 300 & 25 & 1.0942 & 1.67 \\
Pipe 3 (P3) & 500 & 25 & 1.0942 & 1.67 \\ \hline
\end{tabular}
\end{table}

\subsubsection{Computational algorithms and implementation}\label{implementation}

Simulations of \textit{System 1} are performed with two system models, i.e. one based on the FOM as described in \Cref{FOM} and one based on the ROM as described in \Cref{ROM}, that each have a block structure according to \Cref{system1-2} implemented in {\tt Matlab-Simulink} as library blocks. The pipe segments in \Cref{system1-1} are represented by blocks \textit{DHPipe} in \Cref{system1-2} and each determine the output $T_{out}$ as a function of the input $T_{in}$ by either the FOM or ROM for the given system parameters. The blocks in both system models are coupled, i.e. the output $T_{out}$ of block $P1$ forms the input $T_{in}$ for blocks $P2$ and $P3$, meaning that the latter blocks receive an (in principle) continuously varying input. Both the FOM and ROM can handle this. The external input to the system consists of the supply temperature $T_{supply}$ from the producer, defining the input for block $P1$, and the eventual output of the system consists of the output $T_{out}$ of blocks $P2$ and $P3$, defining the (typically different) supply temperatures to user 1 and 2. Note that, since input $T_g$ is constant and identical for all blocks, it does not contribute to the coupling.

The ROM can via input-output relation \eqref{eqTout} describe the response of $T_{out}$ of blocks $P2$ and $P3$ to continuous changes in
inlet temperatures upon using sufficiently short and constant time intervals $\Delta t$ (\Cref{UnitResponse}). However, direct evaluation
via \eqref{eqTout} becomes prohibitively expensive due to the rapidly increasing number of terms in the summations during the progression in time.
A far more efficient way relies on the fact that relation \eqref{eqTout}, upon considering only discrete time levels $t_k = k\Delta t$, admits
reformulation as
\begin{equation}\label{solutionROM2a}
\widetilde{T}_{out}(t_k) = \xvec{F}_1^\dagger \cdot \xvec{\Delta T}_{in} + \xvec{F}_2^\dagger \cdot \xvec{\Delta T}_g,
\end{equation}
with vectors
\begin{equation}\label{solutionROM2b}
\xvec{F}_{1,2} = [F_{1,2}(t_{k})\;\cdots\;F_{1,2}(t_1)]^\dagger,\quad
\xvec{\Delta T}_{in,g} = [\Delta \widetilde{T}_{in,g,0}\;\cdots\;\Delta \widetilde{T}_{in,g,k-1}]^\dagger,
\end{equation}
containing the values of the transfer functions $F_{1,2}(t_i)$ and the temperature changes $(\Delta \widetilde{T}_{in,i},\Delta \widetilde{T}_{g,i})$ as defined in \Cref{UnitResponse}
at the sequences of discrete time levels $t_1\leq t_i\leq t_k$ and $t_0\leq t_i\leq t_{k-1}$, respectively, in (reversed) order. Refer to \ref{appblock} for a detailed derivation. Note that here the trailing term in \eqref{solutionROM2a} simplifies to $F_2(t_k)\widetilde{T}_g$ on account of the uniform and constant $T_g$.

Critical for reliable simulations is an adequate choice for time step $\Delta t$ in \eqref{GlobalFOM} and \eqref{solutionROM2a} for the FOM and ROM, respectively. From a dynamical perspective, two factors must
be considered, namely (i) the characteristic time scale $\tau_{supply}$ of variable supply temperature $T_{supply}$ and (ii) the characteristic time scales
$\tau_{1,2,3}$ of the thermal flows in the pipes. The former and latter require $\Delta t \ll \tau_{supply}$ and $\Delta t \ll \min(\tau_1,\tau_2,\tau_3)$, respectively, for adequately capturing the system
dynamics. The input $T_{supply}$ is constructed from linear interpolation between sampling points separated by time intervals $\Delta t_{res} = 1\;hr = 3600s$ \citep{Jie2012} and thus yields the piece-wise linear profile shown in \Cref{sim1} below with $\tau_{supply} = \Delta t_{res}$. Good estimates for $\tau_{1,2,3}$ are the typical duration of the fastest transient of transfer functions $F_{1,2}$; in \textit{System 1} this occurs for $F_1$ in pipe $P1$ and yields $\min(\tau_1,\tau_2,\tau_3) = \tau_1\approx 50s \ll \tau_{supply}$ (as demonstrated below in \Cref{F1-system1}).

For the FOM -- whether employed for benchmark simulations or identification of the transfer functions $F_{1,2}$ -- both factors
are relevant, since input as well as internal dynamics must be explicitly resolved.
However, for the ROM (given sufficiently accurate identification of $F_{1,2}$ by the FOM) only resolution of the pipe-wise inputs $T_{in,1-3}$ is relevant.
For the ROM of $P1$ this concerns $T_{in,1}=T_{supply}$; for the ROMs of $P2$ and $P3$ this concerns the output of $P1$ (i.e. $T_{in,2,3}=T_{out,1}$) Thus only the first and second factor is relevant in the former and latter case, respectively. This strictly admits component-wise time steps $\Delta t$ within the global ROM of \textit{System 1} yet for simplicity one global $\Delta t$ is adopted and, in consequence, both factors are relevant for the ROM as well.

The above conditions put forth $\Delta t \ll \tau_1 \approx 50 s$ due to the second factor as the most restrictive (and thereby decisive) constraint imposed by the system dynamics.
The FOM is furthermore subject to numerical constraint $\Delta t \leq \Delta t_{max}$ following \eqref{stability}, with here $\Delta t_{max}=2s$.
This yields $\Delta t_{FOM}=2s$ as adequate time step for the FOM that satisfies both constraints \eqref{stability} and $\Delta t \ll 50 s$.
The latter constraint holds also for the ROM and naturally advances $\Delta t_{ROM}=\Delta t_{FOM}=2s$.

\subsubsection{Prediction accuracy of the ROM}\label{predictionaccuracy}

\begin{figure}[h!]
\centering
	\begin{subfigure}[b]{0.45\textwidth}
\includegraphics[width=5.8cm,height=4.2cm]{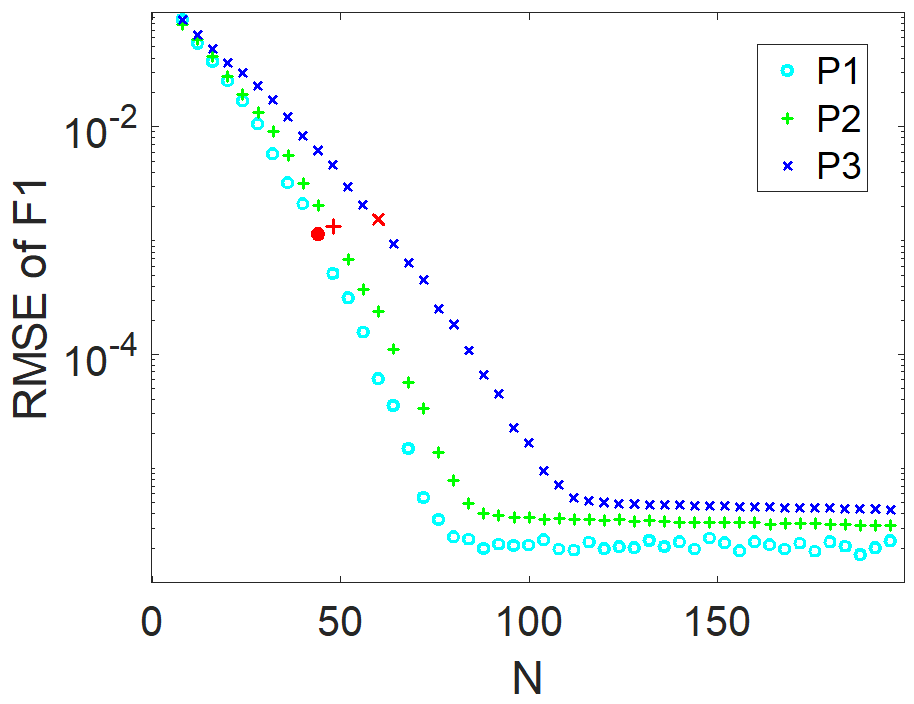}
	\caption{}\label{RMSE-F1}
	\end{subfigure}
	\hfill
	\begin{subfigure}[b]{0.45\textwidth}
\includegraphics[width=5.8cm,height=4.2cm]{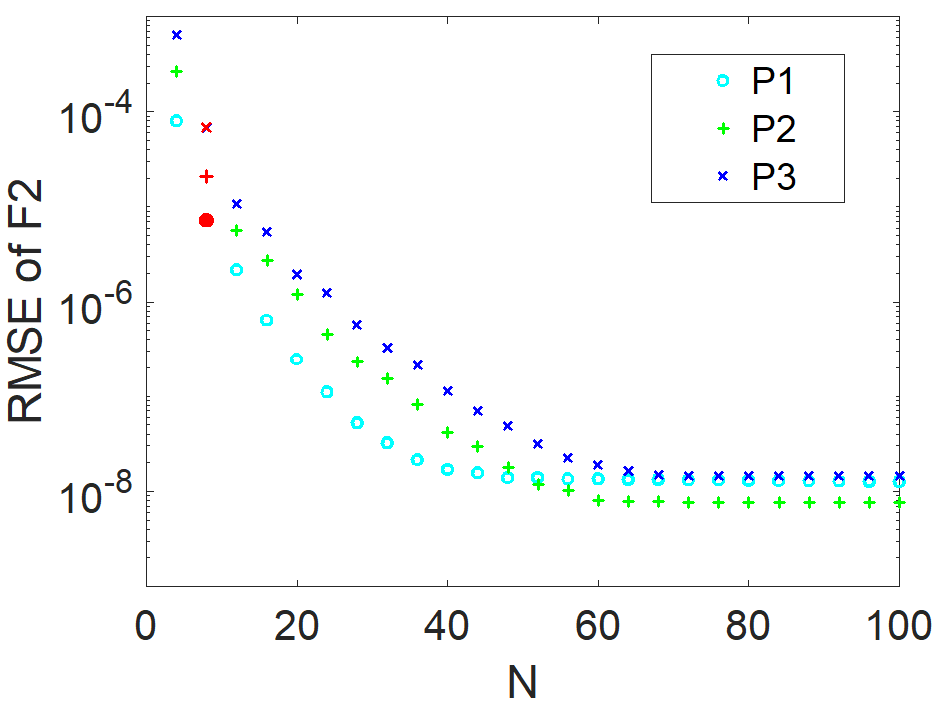}
	\caption{}\label{RMSE-F2}
	\end{subfigure}


\caption{RMSE between the ROM predictions $F_1$ (panel a) and $F_2$ (panel b) and the FOM benchmark for pipes $P1-3$ in \textit{System 1} as a function
of the Chebyshev expansion order $N$. Points highlighted in red correspond with an adequate $N$ for description of $F_1$ and $F_2$ in the ROM.}

\label{RMSE-system1}
\end{figure}

Employment of the ROM starts with identifying transfer functions $F_{1,2}$ for each of the pipes $P1-3$ by the procedure of \Cref{UnitResponse}. Determination of an adequate order $N$ for the corresponding
expansions \eqref{eqCheby} is, as before, done by way of plots of the RMSE \eqref{RMSEdef} similar to \Cref{RMSE}. The RSME profiles for said functions $F_{1,2}$ are shown in \Cref{RMSE-system1} and
obtained from FOM simulations using $\Delta t = 2s$ (\Cref{implementation}). The profiles are qualitatively similar to \Cref{RMSE} yet saturate at substantially higher $N$ compared to the test case
in \Cref{identifications}. This means that $N$ must be increased to reach the same level of prediction accuracy and is a consequence of the pipes in \textit{System 1} being longer than the one in \Cref{identifications}.
\Cref{RMSE-F1} reveals that $F_1$ in $P1$ saturates around $RMSE\approx 2\times{10}^{-6}$ at $N=80$; this onset of saturation shifts to $N=88$ and $N=112$ for $F_1$ in $P2$ and $P3$, respectively, due to
increased pipe lengths (\Cref{tablesystem1}). However, numerical experiments demonstrate that expansion orders yielding $RMSE \approx 1\times{10}^{-3}$ are sufficient for accurate representation of $F_1$
in pipes $P1-3$. This level is already attained at $N =(44,48,60)$ for $(P1,P2,P3)$ (marked in red in \Cref{RMSE-F1}) and these settings are used in the analysis below.

\Cref{RMSE-F2} gives the RMSE profiles for $F_2$ and exposes a dramatic difference with the RMSE profiles of $F_1$ in \Cref{RMSE-F1} by about 2 orders of magnitude. This is comparable with \Cref{RMSE} and
can therefore be attributed to the different response of the system to changes in $T_{in}$ and $T_g$ explained in \Cref{identifications}. Saturation occurs around $RMSE\approx 2\times{10}^{-8}$ yet, similar
as for $F_1$, an RMSE level about 3 order of magnitude above saturation (i.e. $RMSE \approx 1\times{10}^{-5}$) is sufficient for accurate predictions. This yields $N=8$ as adequate expansion order
for $F_2$ in all pipes $P1-3$ (marked in red in \Cref{RMSE-F2}) and also these settings are used hereafter.

\begin{figure}[h!]
\centering
	\begin{subfigure}[b]{0.45\textwidth}
\includegraphics[width=5.7cm,height=4cm]{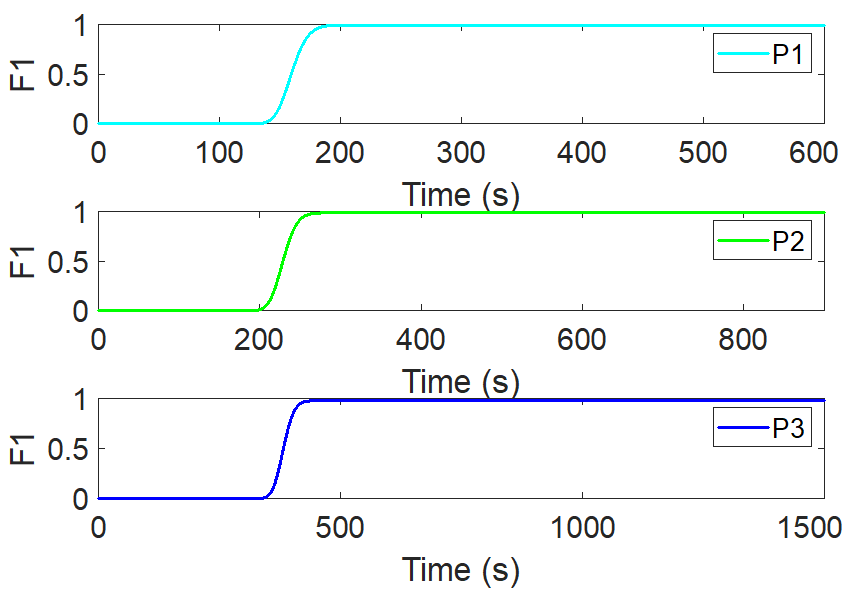}
	\caption{}\label{F1-system1}
	\end{subfigure}
	\hfill
	\begin{subfigure}[b]{0.45\textwidth}
\includegraphics[width=5.7cm,height=4cm]{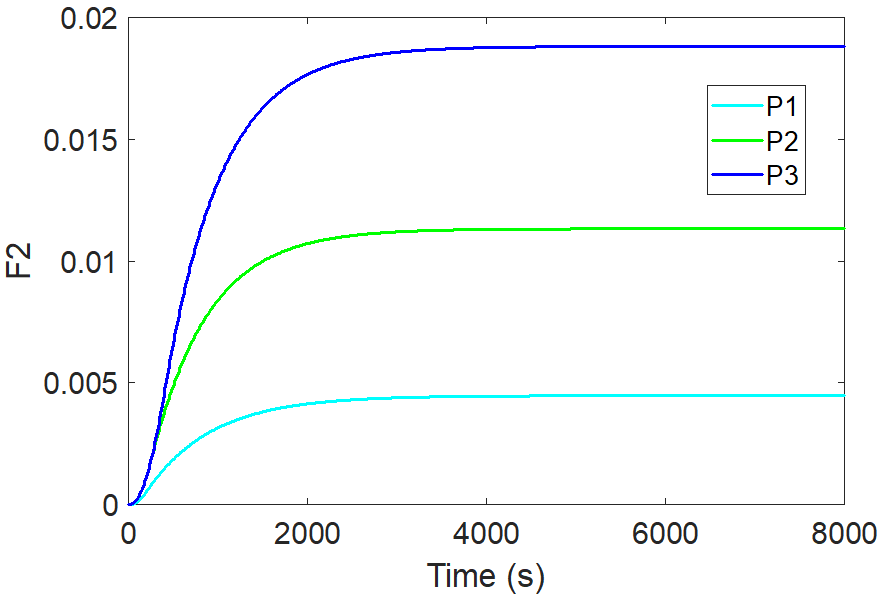}
	\caption{}\label{F2-system1}
	\end{subfigure}
\caption{Identified transfer functions $F_1$ (panel a) and $F_2$ (panel b) for pipes $P1-3$ in \textit{System 1} using expansion orders $N$ as determined from \Cref{RMSE-system1} (highlighted in red).}
\label{F-system1}
\end{figure}

The identified $F_1$ and $F_2$ for the chosen expansion orders $N$ are shown in \Cref{F1-system1} and \Cref{F2-system1}, respectively, and their evolutions resemble those of $F_{1,2}$ in \Cref{fgUnitresponse}.
This further substantiates the above finding that \textit{System 1} exhibits qualitatively the same dynamics as the test case considered in \Cref{identifications}: delayed (yet strong) response to the pipe-wise inlet
temperatures $T_{in,1-3}$ via $F_1$ versus immediate (yet relatively weak) response to ground temperature $T_g$ via $F_2$. Note that the delay increases from $P1$ to $P3$ and thus grows significantly
with pipe length.

The prediction accuracy of the ROM is investigated by comparison with benchmark simulation by the FOM using $\Delta t_{FOM}=2s$ (\Cref{implementation}) for a range of time steps $\Delta t_{ROM}\geq\Delta t_{FOM}$. To this end the departure of the
ROM from the FOM is quantified at time levels $t_k = k\Delta t_{FOM}$ ($k\geq 0$) by RMSE \eqref{RMSEdef}. For mismatching time steps $\Delta t_{ROM}>\Delta t_{FOM}$ the ROM predictions are projected onto time
levels $t_k$ by linear interpolation to always obtain evolutions with a resolution $\Delta t_{FOM}=2s$ as dictated by the system dynamics. This means that (again
given accurate identification of $F_{1,2}$) prediction errors may emanate from insufficient resolution of the input temperatures to pipes $P1-3$ and, in case of said mismatch in time steps, interpolation errors in the output temperatures.

\begin{figure}[h!]
\centering
	\begin{subfigure}[b]{0.45\textwidth}
\includegraphics[width=5.5cm,height=4.3cm]{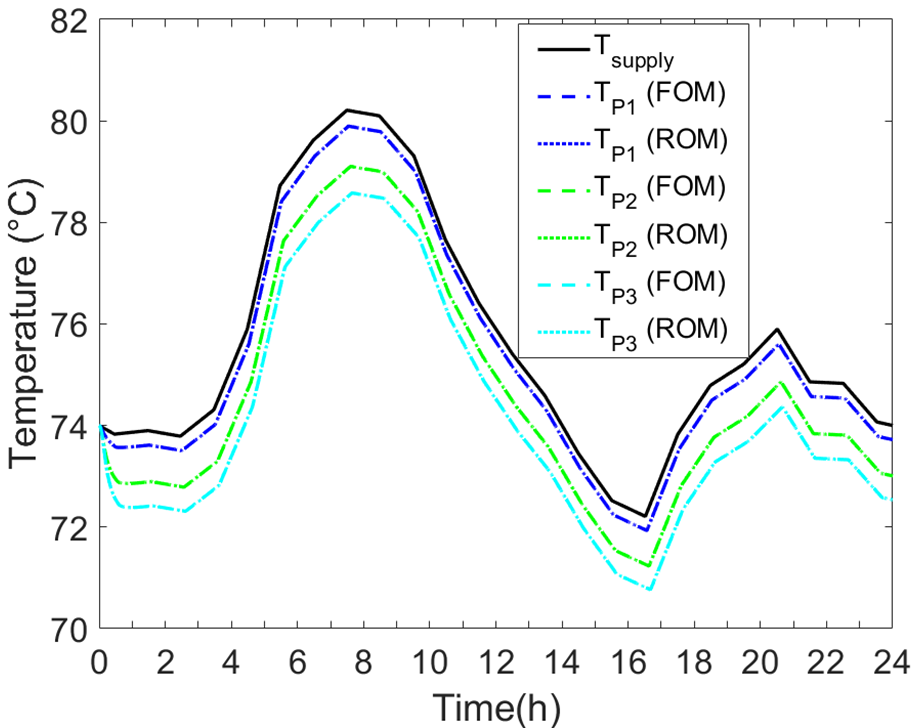}
	\caption{\textit{Sim1 ROM}}\label{sim1}
	\end{subfigure}
	\hfill
	\begin{subfigure}[b]{0.45\textwidth}
\includegraphics[width=5.5cm,height=4.3cm]{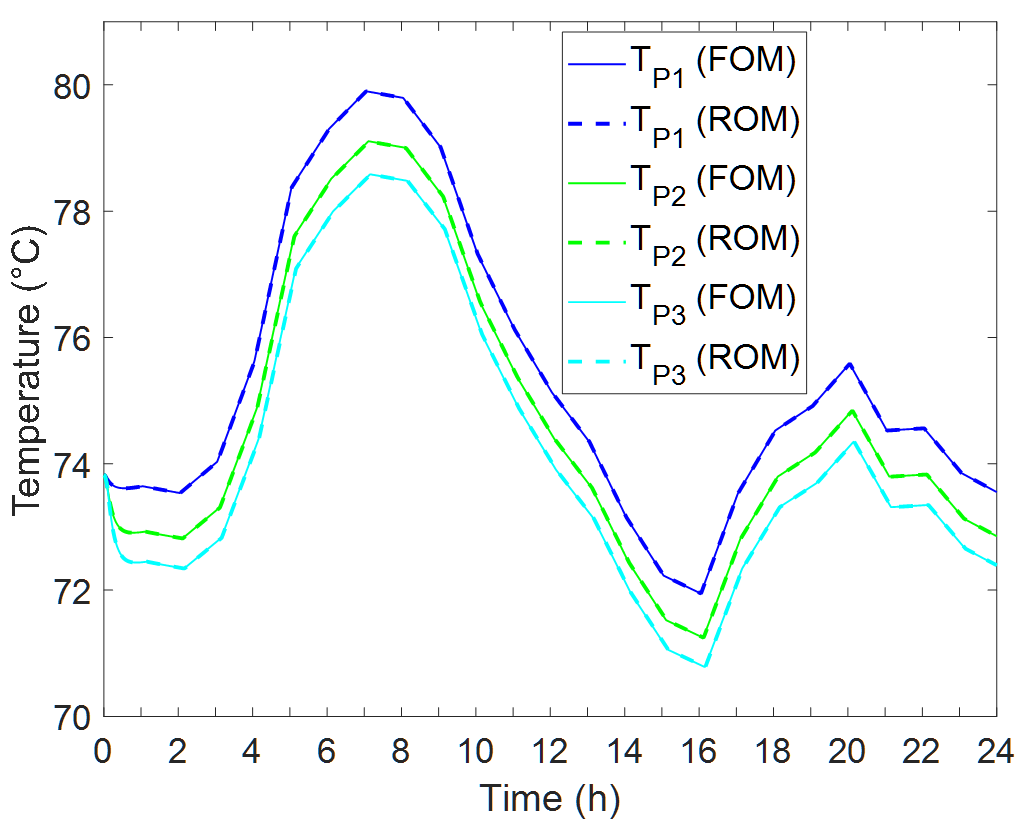}
	\caption{\textit{Sim2 ROM}}\label{sim2}
	\end{subfigure}
	\hfill
	\begin{subfigure}[b]{0.45\textwidth}
\includegraphics[width=5.5cm,height=4.3cm]{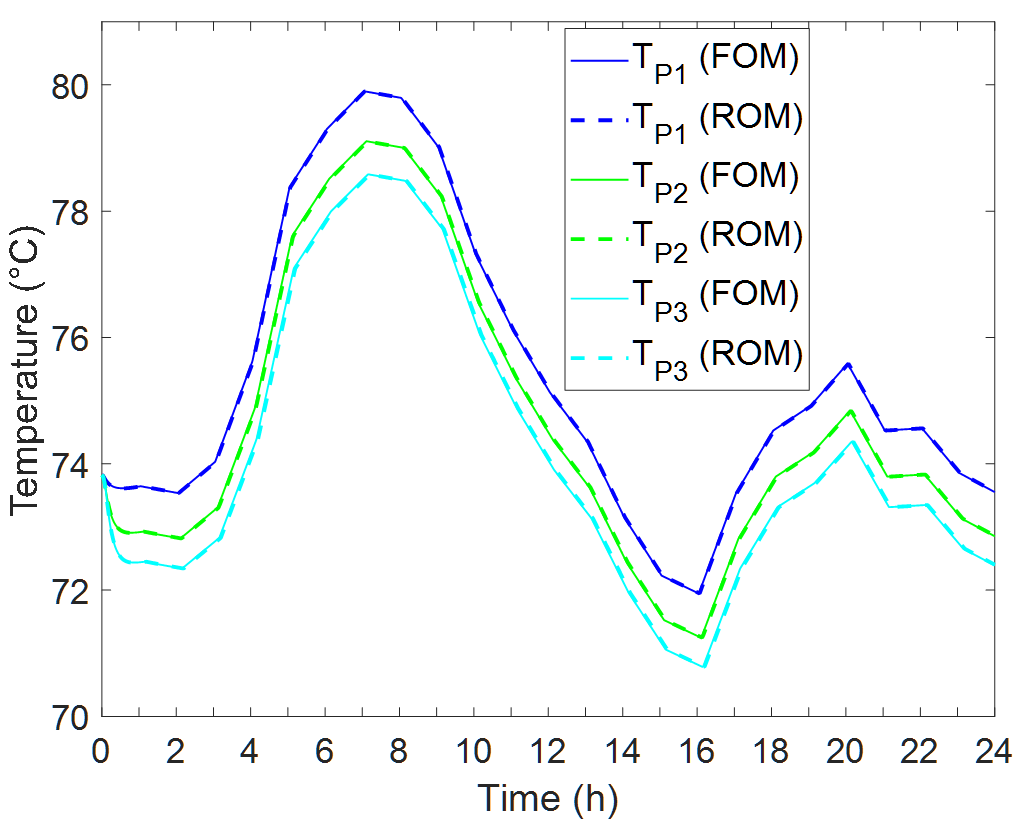}
	\caption{\textit{Sim3 ROM}}\label{sim3}
	\end{subfigure}
	\hfill
	\begin{subfigure}[b]{0.45\textwidth}
\includegraphics[width=5.4cm,height=4.3cm]{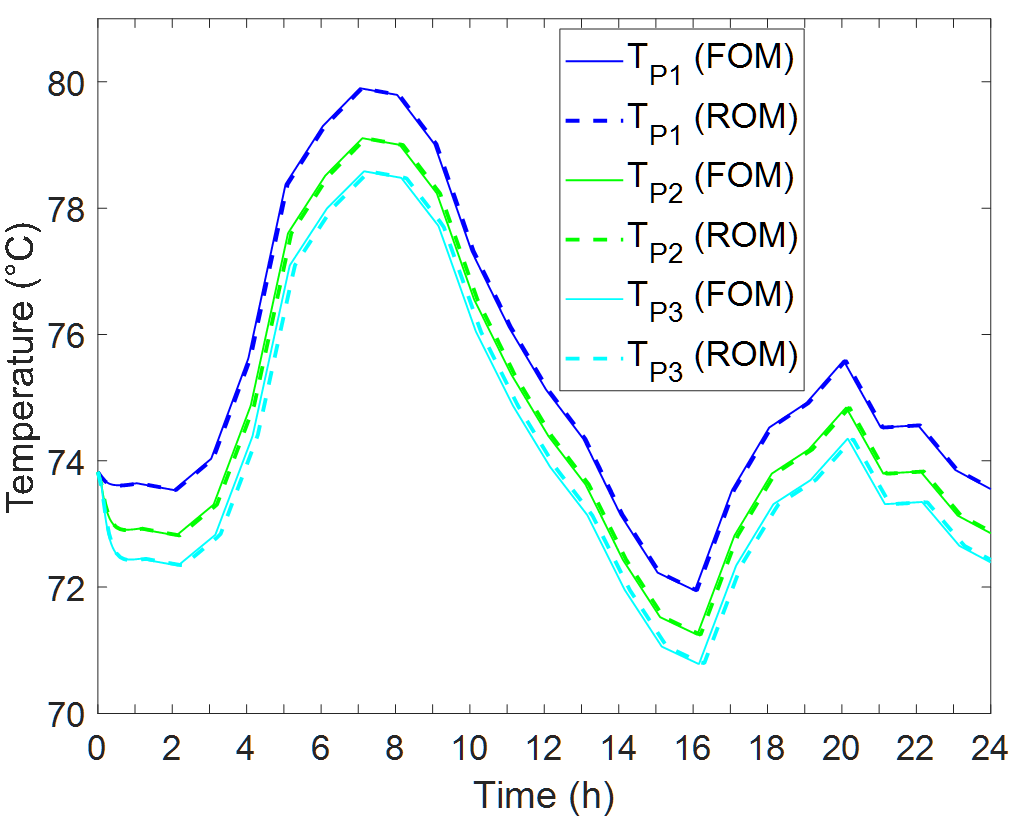}
	\caption{\textit{Sim4 ROM}}\label{sim4}
	\end{subfigure}

\caption{Evolution of the outlet temperatures $T_{P1-3}$ of pipes $P1-3$ in \textit{System 1} for given supply temperature $T_{supply}$ (panel a) according to ROM predictions versus FOM simulations for test cases \textit{Sim1-4 ROM} (\Cref{runs}).}

\label{sims}
\end{figure}

The prediction errors for the individual pipes $P1-3$ are summarized in \Cref{runs} for 4 test cases (\textit{Sim1-4 ROM}) distinguished by time steps ranging from $\Delta t_{ROM}=2s = \Delta t_{FOM}$ to $\Delta t_{ROM}=360s \gg \Delta t_{FOM}$.
This reveals that for $\Delta t_{ROM}=\Delta t_{FOM}$ (\textit{Sim1 ROM}) the prediction error of the ROM is very small for all three pipes and, since output interpolation errors are absent in this case, must be entirely
attributed to finite temporal resolution of the pipe-wise input temperatures. The prediction error is somewhat larger for $P2$ and $P3$ compared to $P1$ and very likely results from the different characteristic
time scales of the respective inputs (i.e. $\tau_{2,3}\sim 50s$ versus $\tau_{supply}=3600s\gg \tau_{2,3}$ for former and latter, respectively). This suggests relatively stronger temporal variation for
$P2$ and $P3$ and, inherently, a greater effect of finite resolution by the ROM. The prediction accuracy nonetheless is high and this is substantiated by the close agreement in \Cref{sim1} between the
output temperatures of the ROM and FOM for given $T_{supply}$ following \citep{Jie2012}.

\begin{table}[]
\centering
\small
\caption{Prediction accuracy of the ROM for pipes $P1-3$ in \textit{System 1} for test cases \textit{Sim1-4 ROM} in terms of the RMSE relative to benchmark \textit{Sim FOM}.}
\label{runs}
\begin{tabular}{lllll}
\hline
Runs & Time step (s) &
  RMSE of P1 &
  RMSE of P2 &
  RMSE of P3 \\\hline
\textit{Sim FOM} & 2 & 0 & 0 & 0\\
\textit{Sim1 ROM} & 2 & $9.94 \times 10^{-4}$ & $3.6 \times 10^{-3}$ & $5.5 \times 10^{-3}$ \\
\textit{Sim2 ROM} & 60 & $5.54 \times 10^{-3}$ & $9.2 \times 10^{-3}$ & $9.7 \times 10^{-3}$ \\
\textit{Sim3 ROM} & 120 & $2.22 \times 10^{-2}$ & $2.95 \times 10^{-2}$ & $4.12 \times 10^{-2}$\\
\textit{Sim4 ROM} & 360 & $5.71 \times 10^{-2}$ & $9.10 \times 10^{-2}$ & $1.37 \times 10^{-1}$\\\hline
\end{tabular}
\end{table}

Test cases \textit{Sim2-4 ROM} reveal that the RMSE increases -- and the prediction accuracy decreases -- with larger time step $\Delta t_{ROM}$. For $P1$ this stems from progressively coarser representation of the piece-wise linear evolution of input $T_{in,1}=T_{supply}$ (\Cref{sim1}) with increasing $\Delta t_{ROM}$ and the projection of the corresponding output $T_{out,1}$ on time levels $t_k$; for $P_2$ and $P_3$ this via $T_{in,2,3}=T_{out,1}$ stems both from the latter and the projection of $T_{out,2,3}$ on $t_k$. The prediction errors following \Cref{runs} grow approximately linear with $\Delta t_{ROM}$ (consistent with employment of linear interpolation schemes) and the prediction error for P2 and P3 also for test cases \textit{Sim2-4 ROM} always exceeds that for $P1$ for reasons explained before. However, for \textit{Sim2 ROM} and \textit{Sim3 ROM} these effects remain virtually invisible in the actual temperature evolutions (\Cref{sim2}-\ref{sim3}); visible departures from the FOM benchmark occur only for \textit{Sim4 ROM} yet nonetheless stay acceptable (\Cref{sim4}).

The above findings imply that the ROM even for coarse time steps of $\Delta t_{ROM}=360s$ can still reliably predict the temperature evolutions within (for practical purposes)
reasonable bounds. This makes the proposed ROM a robust tool for practical simulations of DH systems.


\subsubsection{Computational cost of the ROM}

The above demonstrated the high computational accuracy of the ROM for the simulation of DH systems. A further important aspect is the computational cost of the ROM versus the FOM, which is investigated below in
terms of (i) the floating point operations (flops) that are required to compute the solutions and (ii) the actual runtime of an actual computation.


The FOM propagates the solution per time step $\Delta t$ via the matrix-vector relation \eqref{GlobalFOM} and this involves (per pipe segment) the following operations: matrix-vector multiplication $\xvec{g}_1 = \xvec{S}\xvec{T}_{j}$, with $\xvec{S} = \Delta t \xvec{A} + \xvec{I}$ the pre-computed system matrix, scalar-vector multiplications $\xvec{g}_2 = \Delta t T_{in}\xvec{b}_1$ and $\xvec{g}_3 = \Delta t T_{g}\xvec{b}_2$ and vector additions $\xvec{T}_{j+1} = \xvec{g}_1 + \xvec{g}_2 + \xvec{g}_3$. For a system of $P$ degrees of freedom (DOFs) the associated computational costs are $2P^2$ flops and $P$ flops for matrix-vector and scalar-vector multiplications, respectively, and $P$ flops for vector additions \citep{Hunger2005}. This amounts for \eqref{GlobalFOM} to a total computational cost of $2P^2 + 2P + 2P = 2P(P+2)$ flops per time step and
\begin{equation}\label{costFOM}
C_{FOM} = 2P(P+2)Q \stackrel{P\gg1}{\approx} 2P^2 Q \;\; flops
\end{equation}
for $Q$ time steps, where $P = 4N_x\gg 1$ corresponds with the nodal values of the 4 temperatures.

The ROM propagates the output temperature via relation \eqref{solutionROM2a} by the following operations: vector-vector multiplications $g_1=\xvec{F}_1^\dagger \cdot \xvec{\Delta T}_{in}$ and $g_2=\xvec{F}_2^\dagger \cdot \xvec{\Delta T}_g$, each requiring $2k$ flops for the $k$ DOFs in the intermediate states at time level $t_k$, and scalar addition $\widetilde{T}_{out}(t_k) = g_1+g_2$. This yields a total computational cost of $2k+2k+1 = 4k+1$ flops per time step. However, in general ground temperature $T_g$ can be assumed uniform and constant, simplifying the trailing term in \eqref{solutionROM2a} to the pre-computed scalar $g_2 = F_2(t_k)\widetilde{T}_g$ and reducing the cost to $2k+1$ flops per time step.
The total cost of the ROM for Q time steps thus becomes
\begin{equation}\label{costROM}
C_{ROM} = \sum_{k=1}^Q (2k+1) = \sum_{k=1}^{Q+1} (2k-1)-1 = (Q+1)^2 - 1 \stackrel{Q\gg1}{\approx} (Q+1)^2 \;\; flops,
\end{equation}
and corresponds with an arithmetic sequence of odd numbers \cite{Bronstein2005}.

\begin{figure}[h!]
  \centering
  \begin{minipage}{.55\linewidth}
    \centering
    {\includegraphics[width=6cm,height=7cm]{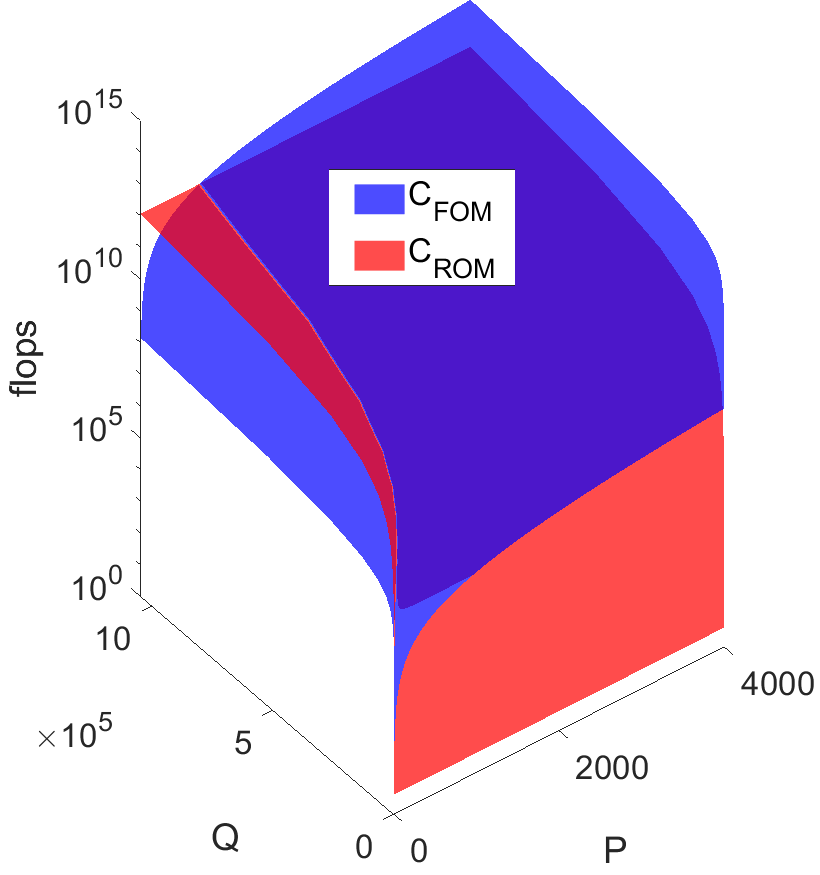}}
    \subcaption{}\label{cost1}
    \end{minipage}\quad
  \begin{minipage}{.42\linewidth}
    \centering
    {\includegraphics[width=\linewidth,height=3.5 cm]{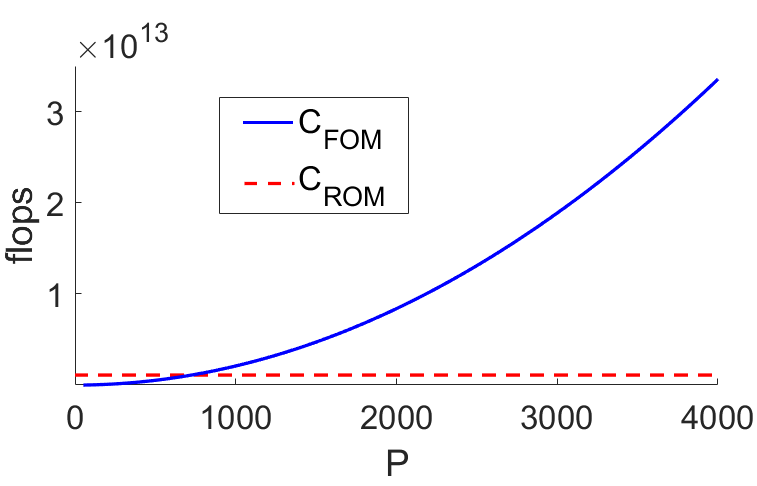}}\subcaption{}\label{cost3}
    {\includegraphics[width=\linewidth,height=3.5cm]{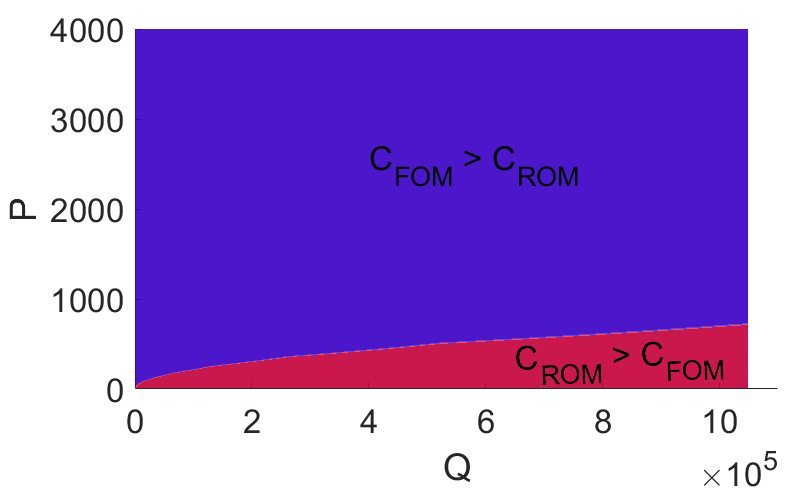}}
    \subcaption{}\label{cost2}
      \end{minipage}
\caption{Computational cost (in flops) for the simulation of one pipe segment by the FOM ($C_{FOM}$) versus the ROM ($C_{ROM}$) as a
function of number of DOFs $P$ and time steps $Q$: (a) $C_{FOM,ROM}$ distribution in $(P,Q)$-space; (b) $C_{FOM,ROM}$ versus $P$ along
cross-section $Q = 10^6$; (c) regions $C_{FOM}>C_{ROM}$ and $C_{ROM}>C_{FOM}$.}
\label{comparisoncost}
\end{figure}

Relations \eqref{costFOM} and \eqref{costROM} reveal that the computational cost of the FOM depends on both spatial resolution $P$ and number of time steps $Q$ while the cost of the ROM is solely dependent on $Q$. \Cref{cost1} gives $C_{FOM}$ and $C_{ROM}$ as a function of $(P,Q)$ in ranges $10\leq P \leq 4000$ and $1\leq Q \leq 10^6$ and clearly demonstrates this essentially different behaviour. The quadratic dependence of
$C_{FOM}$ on $P$, illustrated in \Cref{cost3} for $Q=10^6$, results in a rapid (and monotonic) increase of these costs and causes them to always exceed $C_{ROM}$ beyond the threshold
%
\begin{eqnarray}\label{Pmin}
P_{min}(Q) = (Q+1)/\sqrt{2Q}\approx \sqrt{Q/2},
\end{eqnarray}
i.e. $C_{ROM} > C_{FOM}$ for $P<P_{min}$ and $C_{FOM} > C_{ROM}$ for $P > P_{min}$. This threshold shifts unfavourably for the ROM with longer simulation times due to the quadratic versus linear dependence of $C_{ROM}$ and $C_{FOM}$, respectively, on $Q$. However, \Cref{cost2} demonstrates that region $C_{ROM} > C_{FOM}$ nonetheless remains small compared to region $C_{FOM} > C_{ROM}$ in the considered parameter
regime ($\max P_{min} = 707 \ll \max P = 4000$).
The converse happens in $Q$-direction in that here the ROM outperforms the FOM {\it below} a certain threshold, i.e. $C_{ROM}<C_{FOM}$ for $Q<Q_{max}$ and $C_{FOM}<C_{ROM}$ for $Q>Q_{max}$. This readily follows
from \eqref{costFOM} and \eqref{costROM} and yields
\begin{eqnarray}\label{Qmax}
Q_{max}(P) = P(P+\sqrt{P^2-1}-1)\approx 2P^2,
\end{eqnarray}
as said threshold. However, region $C_{ROM} > C_{FOM}$ is also in $Q$-direction small and in fact non-existent in a substantial part of the parameter regime (i.e. $Q_{max} > \max Q = 10^6$ for $P>707$).
This implies that the ROM is faster than the FOM in regime $P > P_{min}$ and $Q < Q_{max}$ by a factor
\begin{eqnarray}\label{CostRatio1}
\mathcal{R} = \frac{C_{FOM}}{C_{ROM}} = \frac{2P^2Q}{(Q+1)^2}\approx \frac{2P^2}{Q} > 1,
\end{eqnarray}
and may thus structurally reduce the computational costs (e.g. $C = 32$ for $P=4000$ and $Q=10^6$).

\begin{figure}[h]
\centering
\includegraphics[width=7cm, height=4.5cm]{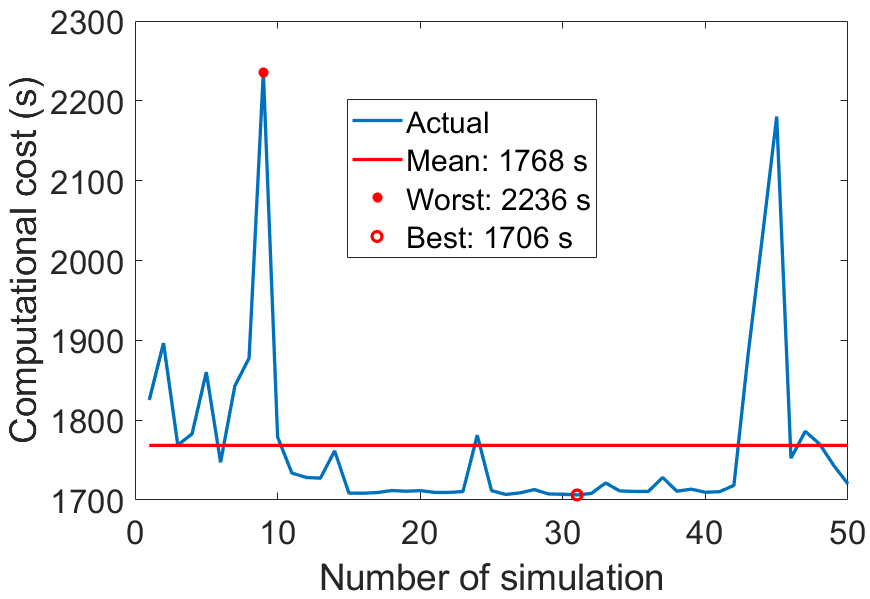}
\caption{Actual runtime $t_{FOM}$ of 50 FOM simulations of \textit{System 1} for $\Delta t = 1 s$ and $\Delta x = 0.1 m$ (yielding $P=40012$ and $Q=86401$) including mean runtime $\bar{t}_{FOM}$ and worst/best cases.}
\label{FOM2s}
\end{figure}

The above cost analysis on the basis of the number of flops associated with the computational schemes \eqref{GlobalFOM} and \eqref{solutionROM2a} includes only the impact of $P$ and/or $Q$. However, the actual runtime of the ROM and FOM (denoted $t_{ROM}$ and $t_{FOM}$, respectively, hereafter) is also dependent on other factors.
The {\tt Matlab-Simulink} simulations with both the ROM and FOM namely involve three stages, i.e. (i) compilation of the {\tt Matlab} source code into block-wise executables, (ii) linking of the executables into one master executable and memory allocation, (iii) actual simulation of the system with the master executable. This process involves more operational, computational and data-handling actions than just flops and makes the runtime slightly different for each simulation, even when re-running a simulation with an already available master executable using the exact same parameter settings. Consider e.g. a FOM simulation of \textit{System 1} for a time span $0\leq t \leq \tau_{sim}$, with $\tau_{sim} = 24 hrs = 86400 s$, using resolutions $\Delta t = 1 s$ and $\Delta x = 0.1 m$. \Cref{FOM2s} gives the actual runtimes $t_{FOM}$ (determined by the {\tt Matlab} commands \textit{tic} and \textit{toc}) for 50 runs on a conventional laptop\footnote{HP ZBook Studio G4 laptop with an Intel core i7-7700 HQ CPU 2.81GHz processor and 32GB installed RAM.} and reveals a significant variation between the worst ($t_{FOM} = 2236 s$) and best ($t_{FOM} = 1706 s$) case around the mean actual runtime $\bar{t}_{FOM} = 1768 s$.
This variation must indeed be attributed to the particular operation and execution of {\tt Matlab-Simulink} elucidated above, since the computational cost \eqref{GlobalFOM} is identical for each run and via $Q = \tau_{sim}/\Delta t +1 = 86401$ and $P = P_1+P_2+P_3 = 40012$ DOFs, amounts to $C_{FOM} = 2.77 \times 10^{14}$ flops. Here the segment-wise DOFs equal $(P_1,P_2,P_3) = (8004,12004,20004)$ and follow from general relations $P = 4N_x$ and $N_x = L/\Delta x + 1$ using the pipe lengths in \Cref{tablesystem1}. Similar variation in actual runtimes occur for {\tt Matlab-Simulink} simulations with the ROM.


The mean actual runtimes $\bar{t}_{ROM}$ and $\bar{t}_{FOM}$ of 50 runs are adopted as estimate for the actual computational cost as a function of $P$ and $Q$. \Cref{costtable1} gives
$\bar{t}_{ROM}$ and $\bar{t}_{FOM}$ versus $P$ for a fixed $\tau_{sim} = 24 hrs = 86400 s$ and $\Delta t = 1s$, which corresponds with a fixed $Q = 86401$.
(Given time step $\Delta t$ is well within the stability bound \eqref{stability}.) Here the computational cost of the FOM increases
from $\bar{t}_{FOM}=25 s$ to $\bar{t}_{FOM}=1768 s$ as the cell size is decreased from $\Delta x=1 m$ to $\Delta x=0.1 m$ and, inherently, the total number of DOFs increases from $P=4012$ to $P=40012$.
The cost of the ROM, on the other hand, remains constant at $\bar{t}_{ROM}=194 s$ due to the beforementioned dependence only on $Q$. This reveals a quadratic dependence
$\bar{t}_{FOM}\propto P^2$ on $P$ for given $Q$ and thus consistency with \eqref{costFOM}. However, $\bar{t}_{ROM}$ and $\bar{t}_{FOM}$ break even at $P_{min}\approx 1.2\times 10^4$ DOFs, which is substantially {\it above} threshold $P_{min} = 208$ DOFs following \eqref{Pmin}. This implies that the performance advantage of the ROM over the FOM is overpredicted by the flops-based analysis. Essentially the same behaviour is found for other $Q$.

\begin{table}[h!]
\centering
\small
\caption{Impact of number of DOFs $P$ on the mean runtime of the FOM ($\bar{t}_{FOM}$) and the ROM ($\bar{t}_{ROM}$) for the simulation of \textit{System 1} for $\Delta t = 1s$ ($Q=86401$).}
\label{costtable1}
\begin{tabular}{ccccccc}
\hline
\multirow{2}{*}{$\Delta x \;(m)$} &
  \multicolumn{3}{c}{Segment-wise DOFs} &
  \multirow{2}{*}{Total DOFs ($P$)} &
  \multirow{2}{*}{$\bar{t}_{FOM}$ (s)} &
  \multirow{2}{*}{$\bar{t}_{ROM}$ (s)}\\ \cline{2-4}
    & \multicolumn{1}{c}{$P_1$}   & \multicolumn{1}{c}{$P_2$}   & $P_3$   &       &       &       \\ \hline
1 & \multicolumn{1}{c}{804} & \multicolumn{1}{c}{1204} & 2004 & 4012 & 25 & 194 \\
0.5 & \multicolumn{1}{c}{1604} & \multicolumn{1}{c}{2404} & 4004 & 8012 & 88 & 194 \\
0.2 & \multicolumn{1}{c}{4004}  & \multicolumn{1}{c}{6004}  & 10004 & 20012 & 424 & 194 \\
0.1 & \multicolumn{1}{c}{8004}  & \multicolumn{1}{c}{12004}  & 20004 & 40012 & 1768 & 194 \\ \hline
\end{tabular}
\end{table}

Investigation of the dependence on the number of time steps $Q$ confirms these observations. \Cref{costtable2} gives $\bar{t}_{ROM}$ and $\bar{t}_{FOM}$ versus $\Delta t$ and $Q$ for a fixed $\tau_{sim} = 24 hrs = 86400 s$ and $\Delta x = 0.5 m$, yielding $P = 8012$ DOFs. (Grid-sensitivity analysis advances $\Delta x = 0.5 m$ as sufficient spatial resolution for \textit{System 1}.) The computational cost of the ROM changes from $\bar{t}_{ROM}=194 s$ to $\bar{t}_{ROM}=0.28 s$ via a quadratic dependence $\bar{t}_{ROM}\propto (Q+1)^2$ on $Q$ for given $P$ and is here thus consistent with \eqref{costROM}. The FOM is examined only for time steps $\Delta t = 1s$ ($Q=86401$) and $\Delta t = 2s$ ($Q=43201$) due to stability criterion \eqref{stability} yet this
nonetheless reveals a proportional dependence on $Q$ that suggests consistency with \eqref{costFOM}. However, break even of $\bar{t}_{ROM}$ and $\bar{t}_{FOM}$ occurs at $Q_{max}\approx 4.5\times 10^4$ time steps,
which is now considerably {\it below} threshold $Q_{max} = 1.28\times 10^8$ following \eqref{Qmax}. This again implies an overprediction of the performance of the ROM by the flops-based analysis. For other $P$ essentially the same behaviour is found.

The above overprediction of the ROM performance results from the overhead caused by operations other than the flops themselves during actual simulations with {\tt Matlab-Simulink}. The actual and flop-based cost namely relate via $\bar{t}_{ROM} = t_{OV} + C_{ROM} f$ for the ROM and (assuming for simplicity comparable overhead) likewise for the FOM, with $t_{OV}$ and $f$ the time consumption by the overhead and per flop, respectively. These relations, upon expressing the overhead as $t_{OV} = C_{OV}f$, yield
\begin{eqnarray}\label{CostRatio2}
\bar{\mathcal{R}} = \frac{\bar{t}_{FOM}}{\bar{t}_{ROM}} = A(\mathcal{R}-1) + 1,\quad\quad A = \frac{C_{ROM}}{C_{OV}+C_{ROM}},
\end{eqnarray}
as the actual improvement factor in terms of the flops-based factor $\mathcal{R}$ according to \eqref{CostRatio1}. Coefficient $A$ decays monotonically from upper bound $\max A = 1$ for $C_{OV}=0$ to asymptotic limit $\lim_{C_{OV}\rightarrow\infty}A=0$ and in conjunction with $\mathcal{R}>1$ implies
\begin{eqnarray}\label{CostRatio3}
\Delta R = \mathcal{R} - \bar{\mathcal{R}} = (1-A)(\mathcal{R}-1) \geq 0.
\end{eqnarray}
thereby demonstrating that a non-zero overhead (i.e. $C_{OV}>0$) in fact {\it always} diminishes the actual cost reduction by the ROM relative to the FOM (i.e. $\Delta R>0$ and
thus $\bar{\mathcal{R}} < \mathcal{R}$). This (at least qualitatively) explains the observed behaviour and exposes minimisation of the overhead as the way to structurally optimize the performance of the ROM.

\begin{table}[h!]
\centering
\small

\caption{Impact of number of time steps $Q$ on the mean runtime of the FOM ($\bar{t}_{FOM}$) and the ROM ($\bar{t}_{ROM}$) for the simulation of \textit{System 1} using $\Delta x = 0.5$ ($P=8012$) for the FOM.}

\label{costtable2}
\begin{tabular}{cccc}\hline
$\Delta t \; (s)$ & Simulation steps (Q) & $\bar{t}_{ROM} (s)$ & $\bar{t}_{FOM} (s)$ \\ \hline
1   & 86401 & 194 & 88 \\
2   & 43201 & 41 & 50 \\
10  & 8641  & 2.3  & - \\
30  & 2881  & 0.6  & - \\
60  & 1441  & 0.42  & - \\
120 & 721   & 0.22  & - \\
360 & 241   & 0.28  & - \\ \hline
\end{tabular}
\end{table}

The above performance analysis concerned identical temporal resolution $\Delta t$ for ROM and FOM. However, a major advantage of the ROM is that $\Delta t$ is determined solely by the required resolution
of the pipe-wise inputs; the FOM on the other hand, is furthermore restricted by numerical stability (\Cref{implementation}). Grid-size sensitivity analysis and stability criterion \eqref{stability}
advance $\Delta x = 0.5 m$ and $\Delta t = 2 s$ as adequate spatio-temporal resolution for the simulation of \textit{System 1} by the FOM. This results in a mean actual runtime $\bar{t}_{FOM} = 50s$
for $\tau_{sim} = 24 hrs$ and corresponding $\bar{t}_{ROM} = 41s$ (\Cref{costtable2}) and, in consequence, an only marginal improvement factor $\bar{\mathcal{R}} = 1.2$ upon
using the same $\Delta t$ for the ROM. \Cref{predictionaccuracy} showed that the ROM enables accurate predictions for time steps as high as $\Delta t = 120 s$, though, and via the corresponding
reduced runtime $\bar{t}_{ROM} = 0.22s$ (\Cref{costtable2}) yields an improvement factor $\bar{\mathcal{R}} = 227$. This signifies a reduction in computational effort by two orders of magnitude
and suggests that the above cost analyses are rather conservative. The actual computational advantage of the ROM may be substantially greater in many practical situations:
\begin{itemize}

\item The FOM in this study concerns the set of 1D systems Eqs.~\eqref{eqwater}-\eqref{eqcasing} and thus involves a spatial discretisation with a relatively small number of
DOFs $P\sim\mathcal{O}(N_x)$. However, this increases dramatically to $P\sim\mathcal{O}(N_x^2)$ and $P\sim\mathcal{O}(N_x^3)$ for DH networks necessitating 2D and 3D pipe
models, respectively, using $N_x$ as typical resolution per coordinate direction. This easily amplifies an improvement factor such as e.g. \eqref{CostRatio1} by several orders of magnitude and effectively
leaves only the ROM as viable option for feasible system simulations.

\item The test cases in the above cost analysis involve time spans $\tau_{sim} = 24hrs = 86400s$. However, in practice often much shorter time spans are needed, such as e.g. hourly predictions in DH system control
 (i.e. $\tau_{sim} = 3600s$). This reduces the number of time steps $Q$ -- and amplifies improvement factor \eqref{CostRatio1} -- by a (further) order of magnitude.

\end{itemize}
The ROM is also much easier to implement and integrate in system models of large DH networks than a FOM of the current (and certainly of a more complex) pipe
configuration. Moreover, any FOM requires a physics-based mathematical model for the system dynamics including detailed information on the geometry and the material properties. Any missing information
makes it very difficult (if not impossible) to build a FOM. The ROM, on the other hand, can be identified purely from data generated by a FOM yet also by (calibration) experiments or field measurements.


\subsection{Design of closed-loop controllers for DH networks}
\label{ControlConsumer}

\begin{figure}[h!]
\centering
	\begin{subfigure}[b]{0.45\textwidth}
\includegraphics[width=6cm,height=2.5cm]{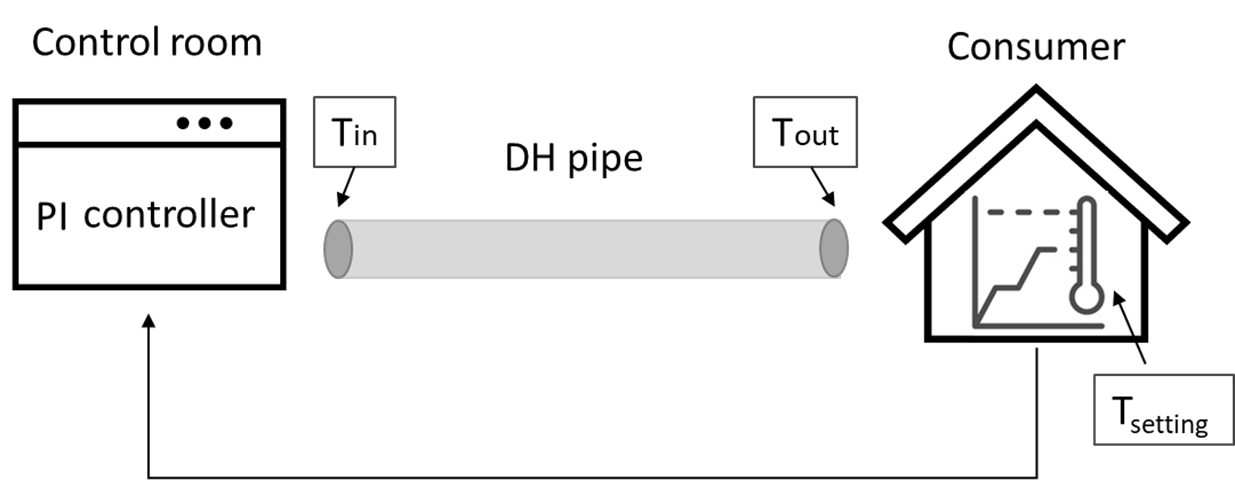}
	\caption{}
	\label{system2-1}
	\end{subfigure}
	\hfill
	\begin{subfigure}[b]{0.45\textwidth}
\includegraphics[width=5.5 cm,height=2.8cm]{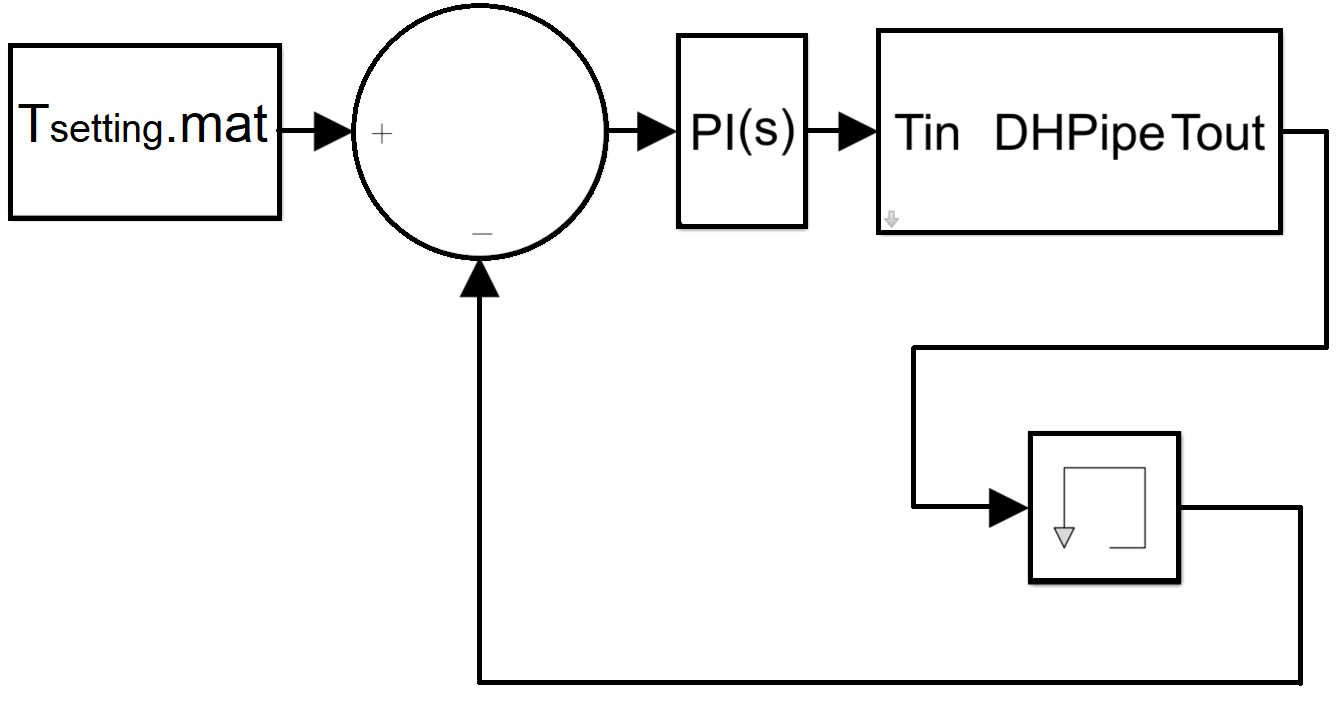}
	\caption{}
	\label{system2-2}
	\end{subfigure}

\caption{User-defined temperature regulation of DH network \textit{System 2}: (a) configuration; (b) {\tt Matlab-Simulink} block structure.}

\label{system2}
\end{figure}

An important application of the ROM is closed-loop control of DH networks. Consider for illustration a simplified system with one energy provider and one user following \Cref{system2-1} (denoted \textit{System 2}) connected by a DN25 pipe of length $L = 100m$ and with mass flux $\dot{m}=1.0942 kg/s$.
The control target is making the outlet temperature $T_{out}$ reach the desired temperature $T_{setting}$ at the consumer side as quickly as possible. This is to be achieved by regulating the inlet temperature $T_{in}$ on the basis of the measured $T_{out}$ via a conventional PI controller as shown in \Cref{system2-2}. Thermal inertia results in a non-trivial response of $T_{out}$ to changes in $T_{in}$ and the ROM, by incorporating the underlying input-output relation (\Cref{ROM}), enables design of a controller (``controller synthesis'') that accomplishes the most effective response to user requirements.
The time step is $\Delta t_{ROM} = 1s$ and thereby deliberately more restrictive than $\Delta t_{ROM} = 2s$ following \Cref{implementation} so as to ensure adequate capturing of the system response to continuous variation in $T_{in}$ due to the control action.

\subsubsection{Control strategy}

The PI controller attempts to minimize the error $e(t) = \widetilde{T}_{setting} - \widetilde{T}_{out}(t)$ between set and actual
temperature by adjusting $\widetilde{T}_{in}$ via the weighted sum of momentary and accumulated error
\begin{equation}
\widetilde{T}_{in}(t) = K_p e(t) + K_i \int_0^t e(\xi)d\xi,
\label{ControlLaw1}
\end{equation}
where $K_p$ and $K_i$ are the weights (denoted ``gains'' in control theory) for the proportional (P) and integral (I) terms, respectively \cite{ControlTheoryEng}.
The control action for a fixed $\widetilde{T}_{setting}$ e.g. establishes the equilibrium $\lim_{t\rightarrow\infty}\widetilde{T}_{out}(t)=\widetilde{T}_\infty = \widetilde{T}_{setting}$ and $\widetilde{T}_{in}$ then converges on $\lim_{t\rightarrow\infty}\widetilde{T}_{in}(t) = \widetilde{T}_{in}^\infty$, where $\widetilde{T}_{in}^\infty$ is defined by the asymptotic state $\widetilde{T}_\infty = \widetilde{T}_{in}^\infty F_1^\infty + \widetilde{T}_gF_2^\infty$ of \eqref{eqrelation}. Reaching $\widetilde{T}_{setting}$ eliminates the {\it momentary} error, i.e. $\lim_{t\rightarrow\infty}e(t)=0$, and thus achieves its sought-after minimisation. The {\it accumulated} error, on the other hand, converges on a non-zero limit value: vanishing of $e(t)$ via \eqref{ControlLaw1} namely implies $\lim_{t\rightarrow\infty}K_i \int_0^t e(\xi)d\xi = \widetilde{T}_{in}^\infty$ and thus
$\lim_{t\rightarrow\infty}\int_0^t e(\xi)d\xi = \widetilde{T}_{in}^\infty/K_i$. Hence its incorporation in \eqref{ControlLaw1} is essential to reach the equilibrium.

The required input $\widetilde{T}_{in}$ at constant mass flux $\dot{m}$ can in practice be achieved in any desired range $\widetilde{T}_{in}^{min}\leq \widetilde{T}_{in}\leq \widetilde{T}_{in}^{max}$ via the mixture of two streams in a small mixing chamber, i.e. stream $\dot{m}_1$ from reservoir 1 at constant $\widetilde{T}_{in}^{min}$ and stream $\dot{m}_2 = \dot{m}-\dot{m}_1$ from reservoir 2 at constant $\widetilde{T}_{in}^{max}$, resulting in
$\widetilde{T}_{in}(t) = (\dot{m}_1(t)\;\widetilde{T}_{in}^{min}+\dot{m}_2(t)\;\widetilde{T}_{in}^{max})/\dot{m}$
as temperature of the mixture that enters \textit{System 2} in \Cref{system2-1}. This approach admits variation of $\widetilde{T}_{in}$ by adjustment of mass fluxes $\dot{m}_{1,2}$ similar to a thermostatic tap, which can be done rapidly and accurately via flow controllers, assuming negligible thermal inertia of the mixing process compared to that of \textit{System 2} itself.

A first step towards optimal control (i.e. reaching the control target by an effective yet also physically achievable and efficient control action) is realized by a basic tuning of the PI controller via the Ziegler-Nichols method \citep{Ziegler1942}. Here the gains in \eqref{ControlLaw1} are determined via
\begin{eqnarray}
K_p = 0.45 K_u,\quad K_i = 0.54 K_u/\tau_u,
\label{ControlLaw2}
\end{eqnarray}
with $K_u$ the ultimate sensitivity when the controller output $\widetilde{T}_{in}(t)$ has stable and consistent oscillations with a period time $\tau_u$ for $K_i=0$. The settings for $K_u$ and $\tau_u$ can be determined heuristically by monitoring $T_{in}$ while changing $K_p$ and this yields $K_u=1.055$ and $\tau_u = 152.7s$ for \textit{System 2}. Through \eqref{ControlLaw2} this gives $K_p = 0.4748$ and $K_i = 0.0037$ as corresponding gains for \eqref{ControlLaw1}.

Set temperature $\widetilde{T}_{setting}$ in its simplest form is a fixed value and the control target then is establishing the (new) equilibrium $\lim_{t\rightarrow\infty}\widetilde{T}_{out}(t)=\widetilde{T}_\infty = T_{setting}$ via control law \eqref{ControlLaw1}.
However, $\widetilde{T}_{setting}$ can also be a prescribed profile in time and the control target then is following this reference as closely as possible (termed ``reference tracking'' in literature \cite{ControlTheoryEng}).
Capability for such reference tracking is relevant for effectively dealing with dynamic heating demands from end-users. Control for both fixed and variable $\widetilde{T}_{setting}$ is considered below in \textit{Scenario 1} and  \textit{Scenario 2}, respectively.

\subsubsection{Design of optimal controllers}

\begin{figure}[t!]
\centering
	\begin{subfigure}[b]{0.45\textwidth}
\includegraphics[width=5.5cm,height=3.7cm]{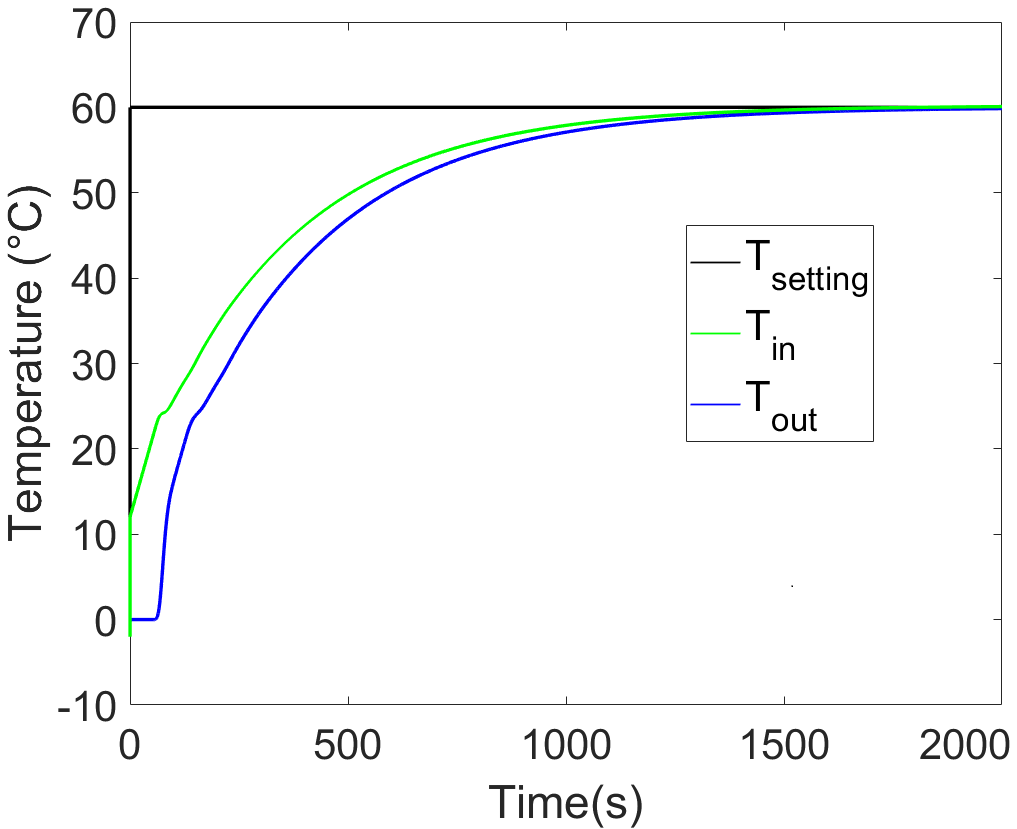}
	\caption{}\label{longresponse}
	\end{subfigure}
	\hfill
	\begin{subfigure}[b]{0.45\textwidth}
\includegraphics[width=5.5cm,height=3.7cm]{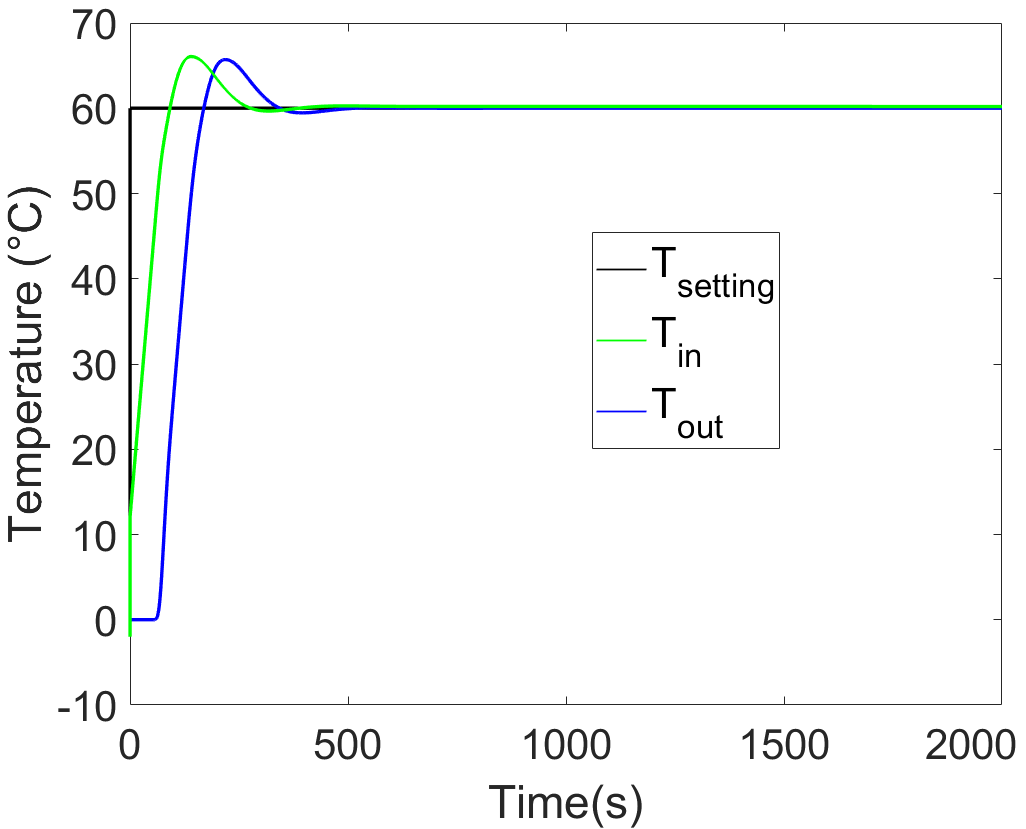}
	\caption{}\label{overshoot}
	\end{subfigure}
\caption{Poor designs of a closed-loop controller for \textit{System 2}: (a) slow response to changes in desired temperature $T_{setting}$ 
for $K_p = 0.211$ and $K_i = 0.002$; (b) overshoot of $T_{setting}$ for $K_p = 0.211$ and $K_i = 0.01$.}
\label{issues}
\end{figure}

\textit{Scenario 1} is examined using $\widetilde{T}_{setting} = 60^\circ C$ as desired temperature and starting from a uniform initial temperature $\widetilde{T}_{out}(0) = \widetilde{T}_0 = 0^\circ C$. The
gains $K_p$ and $K_i$ determined above via the Ziegler-Nichols method are known to be aggressive for many systems in that they may result in strong responses such as e.g. overshoot of the desired temperature $\widetilde{T}_{setting}$ by $\widetilde{T}_{out}$ \citep{Ellis2012}. However, overshoot is often unacceptable for DH systems due to efficiency and safety issues; thermal losses and the risk of (too) hot radiators/convectors both increase with higher $\widetilde{T}_{out}$. This suggests that the basic Ziegler-Nichols tuning is sub-optimal for the current purposes. Its aggressiveness can partially be mitigated by adjusting gain $K_p$ to $K_p = 0.2K_u = 0.211$ as proposed by \citep{McCormack1998} yet optimal control of \textit{System 2} furthermore requires case-specific fine-tuning of gain $K_i$. The latter is rather delicate, though, since increasing $K_i$ decreases the response time -- and thus the effectiveness of the control action -- yet also promotes the unwanted occurrence of overshoot \citep{Ang2005}. This is demonstrated in \Cref{issues}. Gain $K_i = 0.002$ results in an immediate regulation of $\widetilde{T}_{in}$ and, after a short delay due to thermal inertia, in a monotonic evolution of $\widetilde{T}_{out}$ towards $\widetilde{T}_{setting}$ (\Cref{longresponse}). So this controller setting prevents overshoot yet at the expense of a very slow response in that it takes about $1800s = 30min$ for the system to reach $\widetilde{T}_{setting}$. Gain $K_i = 0.01$, on the other hand, triggers a much faster response yet at the cost of a significant overshoot (\Cref{overshoot}).

\begin{figure}[h!]
\centering
	\begin{subfigure}[b]{0.45\textwidth}
\includegraphics[width=5.5cm,height=3.7cm]{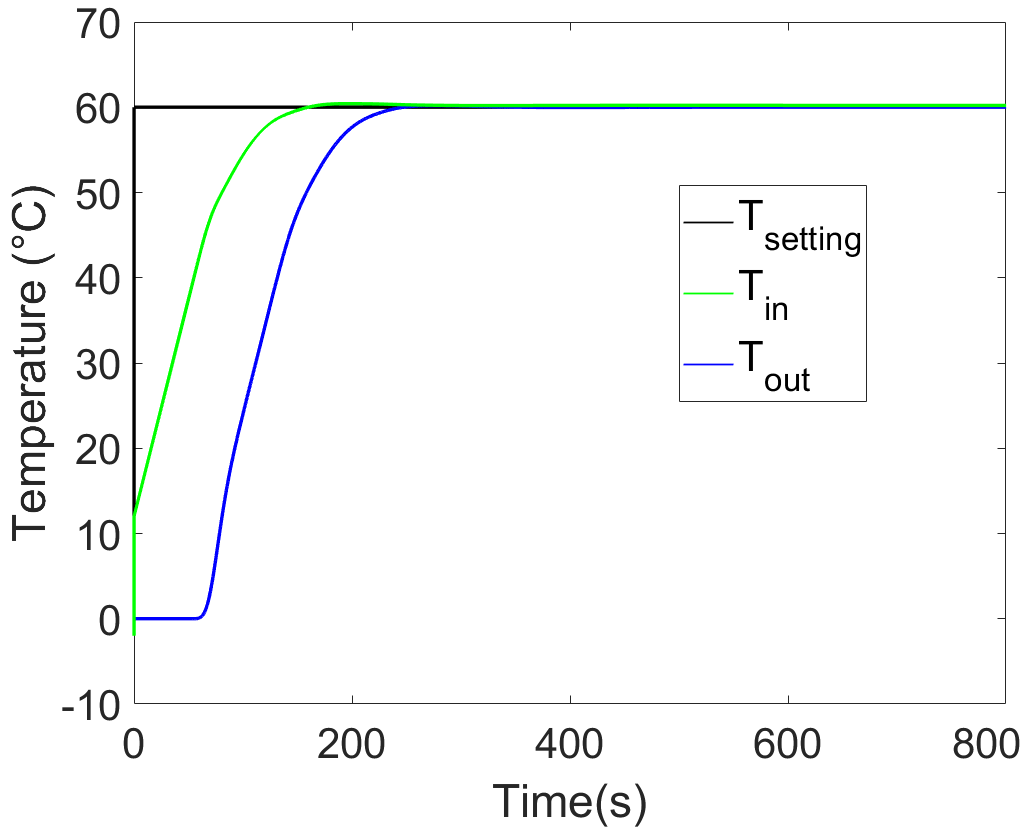}
	\caption{}\label{pid1}
	\end{subfigure}
	\hfill
	\begin{subfigure}[b]{0.45\textwidth}
\includegraphics[width=5.5cm,height=3.7cm]{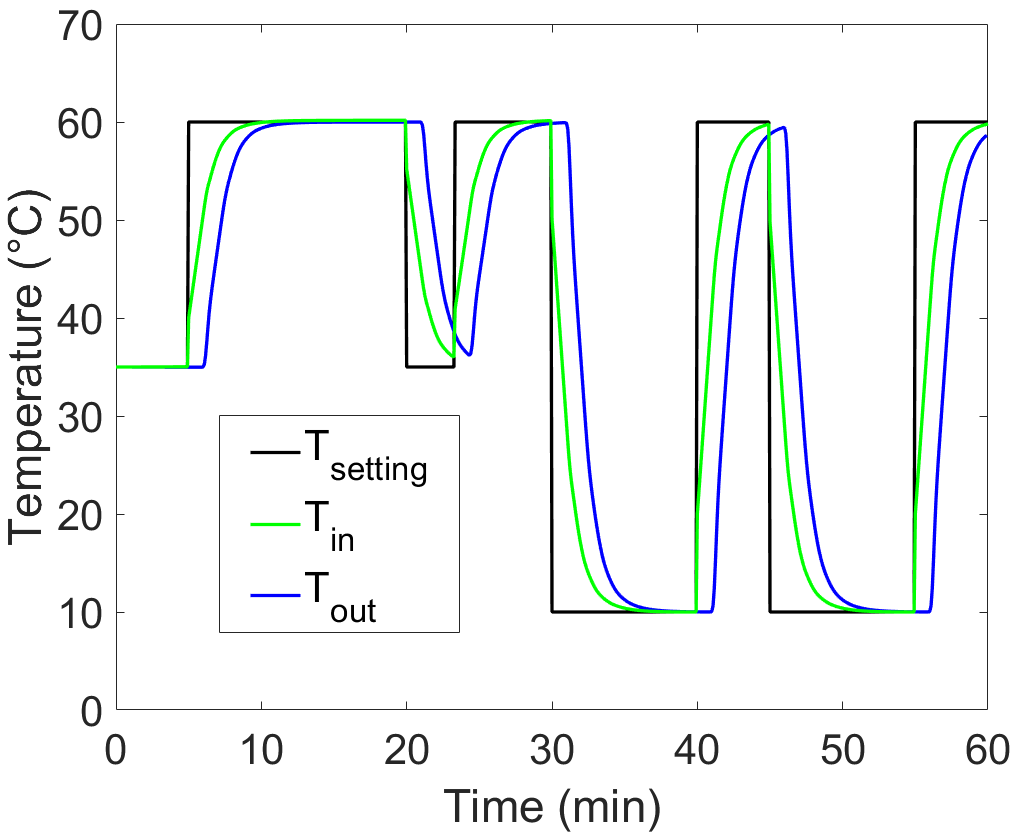}
	\caption{}\label{pid2}
	\end{subfigure}
\caption{Performance of a well-designed controller for \textit{System 2} ($K_p = 0.211$ and $K_i = 0.0085$): (a) response to step-wise change in $T_{setting}$ (\textit{Scenario 1}); (b) reference tracking of 
a user-defined variable $T_{setting}$ (\textit{Scenario 2}).}
\label{pid}
\end{figure}

The behaviour demonstrated in \Cref{issues} suggests that (for given $K_p=0.211$) an optimal gain $K_i$ can be found in the range $0.002\leq K_i \leq 0.01$. The ROM enables efficient and systematic design of an
optimal controller by fine-tuning of $K_i$ via parametric variation. This yields $K_i = 0.0085$ as optimal gain for control law \eqref{ControlLaw1} to accomplish the fastest possible regulation of $\widetilde{T}_{out}$ towards $\widetilde{T}_{setting}$ without overshoot (\Cref{pid1}). Note that the input $\widetilde{T}_{in}$ exhibits some minor overshoot yet this is deemed irrelevant. Striking and in fact counter-intuitive, on the other hand, is that the fine-tuned controller using $K_i = 0.0085$ reaches the desired temperature about twice as fast as the previous ``aggressive'' controller using the higher gain $K_i=0.01$, i.e. $\widetilde{T}_{out}\approx \widetilde{T}_{setting}$ within about $240s = 4 min$ in \Cref{pid1} versus $480s = 8 min$ in \Cref{overshoot}. This underscores the importance and potential of a well-designed controller -- and the need for efficient computational methods such as e.g. the ROM for this purpose -- for an optimal performance of DH networks.

Important to note is that the fine-tuned PI controller aims at {\it optimal} control of \textit{System 2} in the sense mentioned above. However, this is generically not equivalent to
accomplishing the {\it fastest} possible evolution of $\widetilde{T}_{out}$ towards equilibrium $\widetilde{T}_\infty=\widetilde{T}_{setting}$. Consider for illustration the simplified system \eqref{eqcorrelation-simplified0}, which for $h_w=0$ (via $b_2=0$ and $b_1 = -A = \dot{m}/m$) becomes $d\widetilde{T}_{out}/dt = A(\widetilde{T}_{out} - \widetilde{T}_{in}) = b_1(\widetilde{T}_{in} - \widetilde{T}_{out})$.
Analytical solution \eqref{eqcorrelation-simplified} becomes $\widetilde{T}_{out}(t)=U(t)\widetilde{T}_\infty = F_1(t)\widetilde{T}_\infty = F_1(t)\widetilde{T}_{in}$ and corresponds with an instantaneous jump in inlet temperature from $\widetilde{T}_{in}=0$ to $\widetilde{T}_{in}=\widetilde{T}_\infty$ at $t=0$. Smooth transition from $\widetilde{T}_{in}(0)=0$ to $\lim_{t\rightarrow\infty}\widetilde{T}_{in}(t)=\widetilde{T}_\infty$
via a variable input $\widetilde{T}_{in}(t)$, on the other hand, yields an output following
\begin{eqnarray}
\widetilde{T}_{out}(t) = b_1\int_0^t\mbox{e}^{(A(t-\xi)}\widetilde{T}_{in}(\xi)d\xi \quad<\quad U(t)\widetilde{T}_\infty,
\label{eqcorrelation-simplified1}
\end{eqnarray}
for any input without overshoot (i.e. $\widetilde{T}_{in}(t)<\widetilde{T}_\infty$ for any finite $t$) \citep{ControlTheoryEng}. Inequality \eqref{eqcorrelation-simplified1} thus implies that
said jump in $\widetilde{T}_{in}$ establishes the equilibrium faster than the PI controller and in that sense outperforms the latter. However, an instantaneous temperature jump requires
instantaneous jumps in mass fluxes $\dot{m}_{1,2}$ in the beforementioned mixing process for regulating $\widetilde{T}_{in}(t)$. Similarly, reaching or exceeding upper bound $U(t)\widetilde{T}_\infty$ by a variable $\widetilde{T}_{in}(t)$ requires overshoot by the input. Such actions are sub-optimal (if not physically impossible) from a practical perspective and thereby no alternative
for the PI controller. It can be shown that the same reasoning holds for \textit{System 2}.

\textit{Scenario 2} concerns reference tracking and prescribes a variable $\widetilde{T}_{setting}$ for \textit{System 2} to mimic a practical situation of a DH system. Here the heat supply from the DH system is used either for space heating (requiring $\widetilde{T}_{out}=10^\circ C$ and $\widetilde{T}_{out}=35^\circ C$ for moderate and fast heating, respectively) or for domestic hot water (requiring $\widetilde{T}_{out}=60^\circ C$). Set temperature $\widetilde{T}_{setting}$ according to \Cref{pid2} gives a typical hourly demand profile for these purposes and consists of repeated switching between these temperature levels. The PI controller with the fine-tuned gains $K_p=0.211$ and $K_i = 0.0085$ is employed to regulate inlet temperature $\widetilde{T}_{in}$ such that $\widetilde{T}_{out}$ follows this demand profile as closely as possible. The control action $\widetilde{T}_{in}$ and corresponding response $\widetilde{T}_{out}$ is for each step-wise change in $\widetilde{T}_{setting}$ similar to that of the single response in \Cref{pid1} in that $\widetilde{T}_{in}$ is adjusted such that $\widetilde{T}_{out}$ evolves (with the same delay as before) towards the new desired temperature within a time span of about $4mins$ without
over- or undershoot. The controller overall accomplishes a supply temperature $\widetilde{T}_{out}$ that closely follows the demand profile; only the responses to the step-wise changes at $t=20min$ and $t=40min$ somewhat
lag behind in that $\widetilde{T}_{out}$ does not fully reach the desired $\widetilde{T}_{setting}$ within the prescribed interval. This must be attributed to the short duration of this particular temperature
demand compared to the thermal inertia of the system. The ROM in principle enables further fine-tuning of the controller specifically for \textit{Scenario 2} yet this is not done here for
brevity. This case study nonetheless demonstrates that the ROM enables design of optimal controllers also for the (far) more challenging task of reference tracking.

The above employed the ROM for efficient design and fine-tuning of the PI controller. However, an important next step exists in more advanced control strategies based on prediction of future behaviour such as e.g. Model Predictive Control (MPC). MPC is an established and proven approach for process control in fluids and chemical engineering and, given their reliance on similar physical transport phenomena, may be a promising approach also for the control of (sustainable) energy systems such as e.g DH networks \citep{Camacho2013,Corriou2018}. Key for MPC is fast and accurate prediction of
the future system behaviour and the ROM developed here is well-suited for this purpose and may thus enable advanced control of DH systems. This will be addressed in follow-up studies.

\section{Conclusions and future work}\label{conclusion}

This study concerns the development of a data-based compact model for the prediction of the fluid temperature evolution in district heating
(DH) pipeline networks. This consists of a so-called ``reduced-order model'' (ROM) obtained from reduction of the conservation law for energy for each pipe segment to an input-output relation between the pipe outlet temperature and the pipe inlet and ground temperatures. Linearity of the system enables construction of the input-output relation for generic inlet/ground temperature profiles from superposition
of basic input-output relations denoted ``unit step responses'' that can be identified from training data. These step responses are expressed in semi-analytical Chebyshev expansions, which, besides efficient representation, has important advantages as admitting evaluation of function values at arbitrary time levels, easy incorporation in larger system models and exact performance of mathematical operations.

Training data for the ROM are generated by a full-order finite-difference model (FOM) for a 1D pipe configuration. However, the ROM
readily generalizes to more complex pipe configurations involving 3D unsteady heat transfer and 3D steady flow. The only condition for its validity basically is heat-transfer mechanisms being linearly dependent on the temperature field. Thus the ROM in fact holds for a wide range of systems far beyond that considered in the present study.

Performance tests for a single pipe segment reveal that training by FOM data yields a ROM that indeed can accurately describe the input-output relations for arbitrary input profiles. Essential to this end is that the time steps for sampling these input profiles are sufficiently small to adequately capture their temporal variations. The Chebyshev expansions exhibit spectral convergence and thereby demonstrate that this representation is indeed highly suitable for efficiently capturing the input-output relations and constructing a data-based compact ROM.

Performance of the ROM is further investigated for the fast simulation of a small DH network consisting of one pipe segment receiving water with a prescribed supply temperature and a downstream bifurcation into two pipe
segments connected with end users.
Comparison with benchmark simulations by the FOM for a variable supply profile
reveals an accurate prediction of the outlet temperatures of each of the pipe segments by the ROM. Predicted outlet temperatures remain within (for practical purposes) reasonable bounds even for relatively
coarse time steps and this makes the proposed ROM a robust simulation tool for practical DH systems.

The basic computational cost in terms of floating-point operations (flops) of the ROM is less for the FOM only for ``higher'' spatial resolution of the latter and a ``lower'' number of time steps. Moreover, the advantage of the ROM diminishes upon accounting for overhead other than the flops themselves. This suggests that the gain in computational performance by the ROM is limited. However, these findings hold only for relatively simple configurations (e.g. the 1D pipe segments considered here) and in case ROM and FOM require equal temporal resolution.

The situation dramatically changes in favour of the ROM in many practical cases. First, the temporal
resolution for the ROM is determined solely by the pipe-wise input profiles and may therefore often be much coarser than for the FOM (which in addition depends on internal
dynamics and numerical stability). Second, the computational cost for the FOM explodes for spatially more complex systems;
for the ROM this is, regardless of the spatial complexity of the underlying physical system, irrelevant due to its dependence only on time. These factors may quickly amount to a reduction in computational cost by orders of magnitude. Further advantages of the ROM beyond computational aspects exist in ease of implementation and the fact that training data may also be obtained from (calibration) experiments or field measurements.

First successful application of the ROM for control purposes consists of the design and optimization of a PI controller to achieve a user-defined temperature profile at the outlet of a single-pipe
DH system by regulation of the supply temperature. An important next step is incorporation of the ROM in advanced control strategies based on prediction of future behaviour such as e.g. Model Predictive Control (MPC). 
Another important next step (particularly for control purposes) is development of a ROM that, besides inlet/ground temperatures, admits a variable mass flux in pipe segments as a further system input. This fundamentally changes the problem by transforming the pipe configuration into a so-called ``bilinear system'' and, as a consequence, precludes the construction of the input-output relations
from unit step responses. Incorporation of variable mass fluxes thus constitutes a major challenge to further development of the ROM.

\bibliography{Draft5}
\appendix{}
\section{Analytical solution for the semi-discrete model}\label{appsolution}

The evolution of the state vector \eqref{StateVector} with the nodal temperature values is governed by the LTI form \eqref{semi-discrete} of the
energy balances and admits
\begin{equation}\label{tinfsolution}
\xvec{A}\xvec{T}_\infty + \xvec{B}\xvec{g} = 0 \quad\quad\Rightarrow\quad\quad \xvec{T}_\infty = -\xvec{A}^{-1}\xvec{B}\xvec{g} = -\xvec{A}^{-1}(T_{in}\xvec{b}_1 + T_{g}\xvec{b}_2),
\end{equation}
as equilibrium $\lim_{t\rightarrow}\xvec{T}(t) = \xvec{T}_\infty$. Expressing the state vector relative to equilibrium \eqref{tinfsolution}, i.e.
$\xvec{T}' = \xvec{T} - \xvec{T}_\infty$, enables reformulation of \eqref{semi-discrete} as
\begin{equation}
\frac{d\xvec{T}'}{dt} = \xvec{A}\xvec{T}',
\end{equation}
and readily yields
\begin{equation}\label{trelanalytic}
\xvec{T}'(t) = e^{\xvec{A}t}\xvec{T}'_0,
\end{equation}
as analytical solution for $\xvec{T}'$, with $\xvec{T}'_0 = \xvec{T}_0 - \xvec{T}_\infty$ as initial condition. This, in turn, yields
\begin{equation}\label{tt2}
\xvec{T}(t) = \xvec{U}(t)\xvec{T}_\infty + e^{\xvec{A}t}\xvec{T}_0,\quad\quad \xvec{U}(t) = \xvec{I} - e^{\xvec{A}t},
\end{equation}
as analytical solution for the original state vector. Subsequently expressing state vector $\xvec{T}$ relative to its initial condition $\xvec{T}_0$,
i.e. $\widetilde{\xvec{T}} = \xvec{T} - \xvec{T}_0$, leads to
\begin{equation}\label{tt3}
\xvec{\widetilde{T}}(t) = \xvec{U}(t)\widetilde{\xvec{T}}_\infty,
\end{equation}
where
\begin{equation}\label{tt3b}
\xvec{\widetilde{T}}_\infty = \xvec{T}_\infty -\xvec{T}_0 = -\xvec{A}^{-1}(T_{in}\xvec{b}_1 + T_{g}\xvec{b}_2) - \xvec{T}_0 = -\xvec{A}^{-1}(\widetilde{T}_{in}\xvec{b}_1 + \widetilde{T}_{g}\xvec{b}_2),
\end{equation}
follows from expressing \eqref{tinfsolution} also relative to $\xvec{T}_0$. Important to note is that the rightmost form of \eqref{tt3b} only exists upon assuming a uniform initial condition $\xvec{T}_0 = T_0\xvec{1}$. Expansion of \eqref{tt3} via \eqref{tt3b} in terms of $\widetilde{T}_{in,g} = T_{in,g} - T_0$ then gives
\begin{equation}\label{tt4}
\xvec{\widetilde{T}}(t) = -\xvec{U}(t)\xvec{A}^{-1}(\widetilde{T}_{in}\xvec{b}_1 + \widetilde{T}_{g}\xvec{b}_2) = \widetilde{T}_{in}\xvec{f}_1(t) + \widetilde{T}_g \xvec{f}_2(t),
\end{equation}
with $\xvec{f}_1(t) = -\xvec{U}(t)\xvec{A}^{-1}\xvec{b}_1$ and $\xvec{f}_2(t) = -\xvec{U}(t)\xvec{A}^{-1}\xvec{b}_2$.

\section{Implementation of input-output relation in {\tt Matlab-Simulink}}\label{appblock}

Input-output relation \eqref{eqTout} admits efficient implementation in {\tt Matlab-Simulink} library blocks by an expression into a vector form. Represent \eqref{eqTout} to this end by
the generic form
\begin{equation}
T_{out}(t) = T_0 + \sum_{j=0}^{K} H(t-t_j)F(t-t_j)\Delta T_j,
\label{VectorFormROM1}
\end{equation}

with $\Delta T_0 = T_{in,0}$ and $\Delta T_j = T_{in,j} - T_{in,j-1}$ for $j\geq 1$. Output $T_{out}$ is at each discrete time level $t_k = k\Delta t$, with $k\geq 0$ and $t_0=0$,
given by
\begin{equation}
T_{out}(t_k) = T_0 + \sum_{j=0}^{K} H(t_k-t_j)F_1(t_k-t_j)\Delta T_j,
\label{VectorFormROM2}
\end{equation}
and upon using $H(t_k-t_j)=0$ for $j\geq k$ by virtue of \Cref{Heaviside} yields
\begin{eqnarray}
T_{out}(t_0) &=& T_0\nonumber\\
T_{out}(t_1) &=& T_0 + F(t_1-t_0)\Delta T_0 = T_0 + F(\Delta t)\Delta T_0\nonumber\\
T_{out}(t_2) &=& T_0 + F(t_2-t_0)\Delta T_0 + F(t_2-t_1)\Delta T_1 =  T_0 + F(2\Delta t)\Delta T_0 + F(\Delta t)\Delta T_1\nonumber\\
T_{out}(t_3) &=& T_0 + F(3\Delta t)\Delta T_0 + F(2\Delta t)\Delta T_1 + F(\Delta t)\Delta T_2\nonumber\\
&\vdots&\nonumber\\
T_{out}(t_k) &=& T_0 + F(k\Delta t)\Delta T_0 + \cdots + F(j\Delta t)\Delta T_{k-j} + \cdots + F(\Delta t)\Delta T_{k-1},
\label{VectorFormROM3}
\end{eqnarray}
as sequence of outputs for all time levels up to the current time level $t_k$. This leads to
\begin{equation}
T_{out}(t_k) = T_0 + \sum_{j=1}^{k} F(j\Delta t)\Delta T_{k-j} = T_0 + \sum_{j=1}^{k} F(t_j)\Delta T_{k-j}  = T_0 + \xvec{F}^\dagger \cdot \xvec{\Delta T},
\label{VectorFormROM4}
\end{equation}
as expression of \eqref{VectorFormROM2} by the inner product of vector $\xvec{F} = [F(t_k)\;F(t_{k-1})\;\cdots\;F(t_1)]^\dagger$ and
vector $\xvec{\Delta T} = [\Delta T_0\;\Delta T_1\;\cdots\;\Delta T_{k-1}]^\dagger$. Vector $\xvec{F}$ simply contains the values of transfer function $F(t)$ at the discrete time levels up to $t_k$ and can
easily be extracted from the pre-computed total vector $\xvec{F}_{tot} = [F(t_K)\;F(t_{K-1})\cdots F(t_1)]$ containing these values for the entire time span of interest $0\leq t_k \leq t_{K}$
for a given simulation. Computation of $T_{out}(t_k)$ via vector form \eqref{VectorFormROM4} is far more efficient than via summation \eqref{VectorFormROM2} and therefore the former is implemented in {\tt Matlab-Simulink}.

\end{document}